\input ./style/arxiv-general.cfg
\documentclass[aap,MSNbibl,dvips]{arximspdf}
\makeatletter
   \@ifpackageloaded{graphicx}{}{\usepackage{graphicx}}
\makeatother

\doi{10.1214/14-AAP1073}
\volume{25}
\issue{6}
\pubyear{2015}
\firstpage{3251}
\lastpage{3294}
\docsubty{FLA}

\makeatletter
\renewcommand{\Re}{\operatorname{Re}}
\renewcommand{\Im}{\operatorname{Im}}
\renewcommand{\mid}{|}
\newcommand{\rrvert}{\vert}
\newcommand{\rrVert}{\Vert}
\newcommand{\llvert}{\vert}
\newcommand{\llVert}{\Vert}
\newtheorem{theorem}{Theorem}
\newtheorem{proposition}{Proposition}
\newtheorem{lemma}{Lemma}
\newproclaim{ass}{Assumption}
\newproclaim{remark}{Remark}
\newtheorem{cor}{Corollary}
\newcommand{\EE}{\mathbb E}
\newcommand{\RR}{\mathbb R}
\newcommand{\PP}{\mathbb P}
\def\bone{\mathbf{1}}
\makeatother

\begin{document}
\begin{frontmatter}

\title{Weak reflection principle for L\'evy processes}
\runtitle{Weak reflection principle for L\'evy processes}

\begin{aug}
\author[A]{\fnms{Erhan}~\snm{Bayraktar}\corref{}\ead[label=e1]{erhan@umich.edu}\thanksref{T1}}
\and
\author[A]{\fnms{Sergey}~\snm{Nadtochiy}\ead[label=e2]{sergeyn@umich.edu}}
\runauthor{E. Bayraktar and S. Nadtochiy}
\affiliation{University of Michigan}
\address[A]{Department of Mathematics\\
University of Michigan\\
530 Church Street\\
Ann Arbor, Michigan 48109\\
USA\\
\printead{e1}\\
\phantom{E-mail: }\printead*{e2}}
\end{aug}
\thankstext{T1}{Supported by NSF Grant
DMS-09-55463.}

%
\received{\smonth{8} \syear{2013}}
%
\revised{\smonth{7} \syear{2014}}

%
\begin{abstract}
In this paper, we develop a new mathematical technique which allows us
to express the joint distribution
of a Markov process and its running maximum (or minimum) through the
marginal distribution of the process itself. This technique is an
extension of the classical reflection principle for Brownian motion,
and it is obtained by weakening the assumptions of symmetry required
for the classical reflection principle to work. We call this method a
weak reflection principle and show that it provides solutions to many
problems for which the classical reflection principle is typically
used. In addition, unlike the classical reflection principle, the new
method works for a much larger class of stochastic processes which, in
particular, do not possess any strong symmetries. Here, we review the
existing results which establish the weak reflection principle for a
large class of time-homogeneous diffusions on a real line and then
proceed to extend this method to the L\'evy processes with one-sided
jumps (subject to some admissibility conditions). Finally, we
demonstrate the applications of the weak reflection principle in
financial mathematics, computational methods and inverse problems.
\end{abstract}

%
\end{frontmatter}

\setcounter{footnote}{1}
\section{Introduction}

\subsection{Classical reflection principle and its applications}\label{subseSRP}

We start with a brief review of the classical reflection principle for
Brownian motion. Denote by $B^x$ the Brownian motion on a real line,
started form $x$. Given arbitrary levels $U>0$ and $K<U$, we can
compute the joint distribution of $B_t=B^0_t$ and its running maximum
$M_t=\sup_{u\in[0,t]}B_u$ as follows:
%
\begin{eqnarray}\label{eqSRP1}
&& \PP ( B_t \leq K, M_t > U )\nonumber
\\
&&\qquad = \PP ( B_{t-T_U + T_U}
\leq K, T_U < t )
= \PP \bigl( B^U_{t-T_U} \leq K,
T_U < t \bigr)
\nonumber\\[-8pt]\\[-8pt]\nonumber
&&\qquad = \PP \bigl( 2U - B^U_{t-T_U} \leq K,
T_U < t \bigr) = \PP ( B_t \geq2U - K, M_t
> U )
\\
&&\qquad = \PP ( B_t \geq2U - K ),\nonumber
\end{eqnarray}
where $T_U$ is the first hitting time of $U$ by $B$, and $B^U_{s} =
B_{T_U + s}$. The above formula first appeared in the work of Bachelier
\cite{Bacheliergames}, followed by a more rigorous treatment, for
example, by L\'evy \cite{Levyarcsin}. Notice that the above
derivations are based on the following well-known properties of
Brownian motion:
\begin{itemize}
\item\textit{Strong Markov property}: $B_{T_U + t} - B_{T_U}$ is a
standard Brownian motion (started from zero), independent of $\mathcal
{F}_{T_U}$ (where the filtration is generated by $B$).

\item\textit{Continuity}: since the paths of $B$ are continuous,
$B_{T_U}=U$, and in view of the above, $B_{T_U + t}$ is a Brownian
motion started from $U$ and independent of $\mathcal{F}_{T_U}$.

\item\textit{Symmetry}: the distribution of $B^U_t$ is symmetric with
respect to the initial level~$U$, that is, $\operatorname{Law}
(B^U_t )=\operatorname{Law} (2U-B^U_t )$.
\end{itemize}

We will come back to the above observations in the next subsection, but
first let us outline several applications of the classical reflection
principle. One obvious application is the computation of the joint
distribution. Since the marginal distribution of a Brownian motion is
available in closed form, the above formula gives us a closed form
expression for the joint distribution of the process and its running
maximum, at any given time.
A more subtle application, which requires the use of reflection
principle itself (rather than the resulting formula for the joint
distribution) comes from financial mathematics. Namely, the reflection
principle turns out to be very useful in the problem of hedging barrier
options. To simplify the notation, throughout the rest of the paper, we
consider $U=0$. Let us assume that the risk-neutral evolution of the
underlying is described by a Brownian motion $B^x$ started from $x\leq
0$ (we assume no discounting). Consider an up-and-out option, written
on this underlying, with a terminal payoff function $h$, such that
$\operatorname{supp}(h) \subset(-\infty,0)$. The payoff of such
barrier option, at the time of maturity $T$, is given by
\[
h\bigl(B^x_T\bigr) \bone_{ \{ \sup_{u\in[0,T]}B^x_u < 0  \}}. %
\]
To find the price of the option, we need to compute the expectation of
the above random variable. This problem can be solved by applying the
formula for the joint distribution of $(B_T,\sup_{u\in[0,T]}B_u)$,
given in (\ref{eqSRP1}). However, making use of the reflection
principle itself, we can obtain more than just a price: in fact, we can
find a static hedging strategy for a given barrier option via the
European-type ones. Recall that a European-type option pays, at the
time of maturity, a certain function of the terminal value of the
underlying. Hence, we need to find a function $G$, such that up until
hitting the barrier, the price of the target barrier option coincides
with the price of a European-type option with maturity $T$ and payoff
$G(B^x_T)$,
\[
\EE \bigl( h\bigl(B^x_T\bigr) \bone_{ \{ \sup_{u\in[0,T]}B^x_u <
0  \}}
\mid \mathcal{F}_{t\wedge T_{0}} \bigr) = \EE \bigl( G\bigl(B^x_T
\bigr) \mid \mathcal{F}_{t\wedge T_{0}} \bigr), %
\]
where the filtration is generated by $B^x$, and $T_0$ is the first
hitting time of zero by~$B^x$.
In fact, due to the tower property, it is enough to ensure that the
above identity holds for $t=T$. To do this, we consider two cases:
$T_{0} \geq T$ and $T_{0} < T$. In the first case, the above equation
reduces to
\[
h\bigl(B^x_T\bigr) \bone_{ \{B^x_T < 0  \}} = G
\bigl(B^x_T\bigr) \bone_{
\{B^x_T < 0  \}}, %
\]
which yields that, below the barrier, the payoff function $G$ has to
coincide with $h$:
$G(z) \bone_{ \{z<0  \}} = h(z)$. Thus we can search for
$G$ in the form $G(z) = h(z) - g(z)$, where $\operatorname
{supp}(g)\subset(0,\infty)$. Considering the second case, and making
use of the \textit{strong Markov property} and \textit{continuity} of $B$, we obtain
%
\begin{equation}
\label{eqSRP2} \bigl(\EE h \bigl(B^0_{\tau} \bigr) - \EE g
\bigl(B^0_{\tau} \bigr) \bigr)_{\tau=T-T\wedge T_0}
\bone_{ \{T_0<T \}} = 0,
\end{equation}
where $B^0_s = B^x_{T_0+s}$ is a new Brownian motion started from zero,
which is independent of $\mathcal{F}_{T_0}$. We emphasize that the
main difficulty in finding $g$ that satisfies~(\ref{eqSRP2}) is the
requirement that \emph{$h$ and $g$ have supports on the opposite sides
of the barrier $U=0$}. In the case of a Brownian motion, we can make
use of its \textit{symmetry}, to conclude that the function $g(z) = h(-z)$
fulfills (\ref{eqSRP2}). In fact, such choice of $g$ satisfies
%
\begin{equation}
\label{eqMIBM} \EE \bigl( h \bigl(B^0_{t} \bigr) \bigr) =
\EE \bigl( g \bigl(B^0_{t} \bigr) \bigr)\qquad\mbox{for
all }t>0,
\end{equation}
which is sufficient for (\ref{eqSRP2}) to hold.
Thus $G(z) = h(z) - h(-z)$ is the solution to the static hedging problem.

It is worth mentioning that the solution to the integral equation (\ref
{eqMIBM}) has an interpretation through the partial differential
equations (PDEs).
Namely, a function $g$ that satisfies (\ref{eqMIBM}) and has support
in $[0,\infty)$ provides a solution to the following inverse problem:
assume that the function $h$, with support in $(-\infty,0]$, serves as
the initial condition to the following heat equation:
%
\begin{equation}
\label{eqinitValBMdef} \cases{ \partial_t u - \Delta u = 0, &\quad $x\in\RR$,
$t>0$,
\cr
u(x,0) = h(x).}
\end{equation}
The Feynman--Kac formula implies that $u(x,t) = \EE h (B^x_t
)$. Then it follows from~(\ref{eqMIBM}) that replacing $h$ with $g$,
we obtain a solution to (\ref{eqinitValBMdef}) which coincides with
the original solution at $x=0$, for all $t>0$. Next, assume that the
initial condition in (\ref{eqinitValBMdef}) has an arbitrary
support in $\RR$. The above observation implies that we can modify
this initial condition (only) in $[0,\infty)$, to ensure that the
resulting solution is zero at $x=0$, for all $t>0$.
One can think of the initial temperature distribution, which we only
control on the positive half line, and only at the initial time. The
goal is then to ensure that the temperature at $x=0$ remains zero at
all times. Of course, in the case of a Brownian motion, the solution to
this problem is obvious due to the symmetry of the system,
$g(x)=h(-x)$. Equivalently, this follows from the symmetry of the
associated heat equation: its first order coefficient is zero, and the
second order coefficient is constant, which corresponds to constant
heat conductivity at all points of the real line.
However, a similar problem can formulated for diffusion processes that
do not posses any symmetries. In this case, the Laplace operator in
(\ref{eqinitValBMdef}) is replaced by a more general elliptic
operator, whose coefficients may not possess the desired symmetries
(e.g., this corresponds to nonconstant and asymmetric heat
conductivity). Although the straightforward approach does not work in
this case, the \emph{weak symmetry mapping}, introduced in the next
section, allows us to solve this problem.

\subsection{Weak reflection principle}\label{subseWRP}

From the above description of the classical reflection principle, it is
easy to see that this method continues to work if we substitute the
assumptions of strong Markov property, continuity and symmetry with the
following: for any $t>0$, the conditional distribution of $B_{T_U +
t}$, given $\mathcal{F}_{T_U}$, is symmetric with respect to the
initial level $U$. Thus the key property is the \textit{symmetry} of the
underlying stochastic process. The strong Markov property and
continuity are only needed to pass from the conditional distributions
to the unconditional ones. This observation immediately yields several
extensions of the classical reflection principle. First, it is well
known that the reflection principle also works for an exponential of a
Brownian motion. Indeed, for any monotone continuous function $F$,
defined on a real line, the process $F(B^x)$ is still strongly Markov
and continuous. Assuming the range of $F$ contains $U=0$, we conclude
that the random variable $X_t=F (B^{V}_t )$, with $V =
F^{-1}(0)$, has the same distribution as $F (2V -
F^{-1}(X_t) )=F (2V - B^V_t )$. Moreover, the function
$x\mapsto F (2V - F^{-1}(x) )$ maps $(-\infty,0]$ into
$[0,\infty)$, and vice versa. Then, for any function $h$, with support
in the negative half line, the function
\[
g(x) = h \bigl(F \bigl(2V - F^{-1}(x) \bigr) \bigr) %
\]
has support in the positive half line and solves the desired integral equation
\[
\EE h (X_t ) =\EE h \bigl(F\bigl(B^V_{t}
\bigr) \bigr) = \EE h \bigl( F \bigl(2V - B^V_t \bigr)
\bigr) = \EE g (X_t )\qquad\mbox{for all }t>0. %
\]
In general, we say that a real-valued stochastic process $X$, with
$X_0=0$, possesses a \textit{strong symmetry} (with respect to zero) if
there exists a mapping $\mathbf{S}\dvtx \RR\rightarrow\RR$, such that
$x\mathbf{S}(x)\leq0$, for all $x\in\RR$, and $\operatorname
{Law}(\mathbf{S}(X_t)) = \operatorname{Law}(X_t)$, for all $t>0$.
Then, for any given function $h$, with support on one side of zero, the
target equation
\[
\EE h (X_t ) = \EE g (X_t )\qquad\mbox{for all } t>0,
\]
has the solution $g(x) = h (\mathbf{S}(x) )$, and the
function $g$ has support on the opposite side of zero. Thus the
reflection principle can be extended easily to any strong Markov
process that possesses a strong symmetry and does not jump across the
barrier $U=0$. In particular, it holds for any diffusion martingale
whose coefficient is an even function, since any such process possesses
a strong symmetry, with $\mathbf{S}\dvtx  x\mapsto-x$.

However, the assumption of the existence of a strong symmetry excludes
many processes important for applications (some of these examples are
discussed in \cite{CarrN}).
In the present paper, we develop a weak formulation of the reflection
principle, which can be applied to a large class of stochastic
processes that \textit{do not posses any strong symmetries}. This new
formulation, albeit weaker than the standard one, is sufficient to
solve the problems outlined in Section~\ref{subseSRP}.
Herein, we restrict our analysis to the strong Markov processes which
do not jump across the given upper (lower) barrier from below (above).
Then, in view of the above discussion, it suffices to consider the
unconditional, as opposed to conditional, distributions of the process.
Consider a stochastic process $X$, defined on a real line and started
from zero. Consider two spaces, $\mathcal{B}^-$ and $\mathcal{B}^+$,
consisting of all Lebesgue measurable functions $h$, such that $h(X_t)$
has finite expectation, for all $t\geq0$, and if $h\in\mathcal
{B}^-$, then $h$ has support in $(-\infty,0]$, while if $h\in\mathcal
{B}^+$, then $h$ has support in $[0,\infty)$. We say that $X$
possesses an \textit{upper weak symmetry} (with respect to zero)\footnote
{Similarly, one can define the weak symmetry with respect to any level
$U$.} if there exists a mapping $\mathbf{W}^+$ from a space of test
functions $\mathcal{B}^-_0\subset\mathcal{B}^-$ into $\mathcal
{B}^+$, such that
\[
\EE h (X_t ) = \EE \bigl(\mathbf{W}^+ h (X_t ) \bigr)
\qquad\mbox{for all } t>0, %
\]
for any $h\in\mathcal{B}^-_0$. Analogously, one can define the \textit{lower weak symmetry}, along with the mapping $\mathbf{W}^-\dvtx \mathcal
{B}^+_0\rightarrow\mathcal{B}^-$. We do not insist on a particular
choice of the spaces of test functions $\mathcal{B}^{\pm}_0$, but, in
what follows, we choose the spaces that include all smooth functions
with compact support. We will refer to $\mathbf{W}^{\pm}$ as the \textit{weak symmetry mapping}, although one should remember that in the
present setting, the specific form of the symmetry may vary
significantly depending on the underlying process~$X$. For example, if
the distribution of $X_t$ is not symmetric with respect to zero, the
image of a piece-wise linear function may be a curve with nonzero
curvature at every point.
Thus the \textit{weak reflection principle} consists of the application of
the strong Markov property, the (semi-) continuity and the weak
symmetry mapping~$\mathbf{W}^{\pm}$.

It is worth mentioning that there is a particular choice of weak
symmetry that has received a lot of attention in the existing
literature.\footnote{We thank the anonymous referee for pointing this
out.} This symmetry is related to the geometric Brownian motion with drift.
Notice that a Brownian motion with a constant (nonzero) drift does not
possess any strong symmetries, and neither does its exponential---the
geometric Brownian motion---which is the main component of the
celebrated Black--Scholes--Merton model.
Nevertheless, Carr and Chou \cite{CarrChou} eliminate the drift by a
Girsanov change of measure and derive the following relation:
%
\begin{equation}
\label{eqspecialWSdef} \EE \biggl( f \biggl(\frac{S_T}{S_{T\wedge\tau}} \biggr) \Big|
\mathcal{F}_{T\wedge\tau} \biggr) = \EE \biggl( \biggl(\frac{S_T}{S_{T\wedge\tau}}
\biggr)^{\alpha} f \biggl(\frac{S_{T\wedge\tau}}{S_T} \biggr) \Big|
\mathcal{F}_{T\wedge\tau} \biggr),
\end{equation}
which holds for any (admissible) function $f$, with $S$ being a
geometric Brownian motion and $\tau$ being the first hitting time of a
(strictly positive) barrier. The value of $\alpha$ is given explicitly
via the drift and volatility of $S$. Relation (\ref{eqspecialWSdef}) is extended to an arbitrary stopping time $\tau$
in \cite{CarrPCS}, where it is also shown that (\ref
{eqspecialWSdef}) still holds if $S$ is given by a geometric
Brownian motion run on an independent continuous stochastic clock. The
equivalent formulations of (\ref{eqspecialWSdef}), and the
additional properties of stochastic processes $S$ that satisfy this
relation are established in \cite{Tehranchi}. Finally, the authors of
\cite{Molchanov} and \cite{Rheinlander} further develop the analysis
of (\ref{eqspecialWSdef}), which they call the (\emph{quasi})
\textit{self-duality}, in the case when $S$ is a L\'evy process and in the
multivariate case.
It is important to observe that, in fact, (\ref{eqspecialWSdef})
specifies a weak symmetry mapping associated with $X=\log(S)$. To see
this, assume that $X$ is a strongly Markov process, which does not jump
across the barrier from below (resp., above), and that $\tau$ is the
first time when $S=\exp(X)$ hits a given upper (resp., lower) barrier.
Without loss of generality, we can also assume that the barrier is
equal to one and that $f$ has support in $(0,1)$ [resp., $(1,\infty)$].
Then, repeating the derivation of (\ref{eqSRP2}), we conclude that
\[
\EE f \bigl(\exp(X_t) \bigr) = \EE \bigl( e^{\alpha X_t} f \bigl(
\exp(-X_t) \bigr) \bigr)\qquad\mbox{for all } t>0, %
\]
is a sufficient condition for (\ref{eqspecialWSdef}) to hold, which
also becomes necessary if the distribution of $\tau$ has full support
in $(0,\infty)$. Thus, in a Markovian setting, (\ref
{eqspecialWSdef}) can be viewed as the following weak symmetry mapping:
%
\begin{equation}
\label{eqspecialWSdef2} \mathbf{W}^{\pm}h(x) = e^{\alpha x} h(-x).
\end{equation}
It is shown in \cite{CarrPCS,Tehranchi,Molchanov}
and \cite{Rheinlander} that many popular stochastic processes admit
(\ref{eqspecialWSdef2}) as the weak symmetry mapping. However, it
is also easy to see (cf. \cite{CarrN}) that there are many important
Markov processes whose weak symmetry mapping is different from (\ref
{eqspecialWSdef2}) and cannot be obtained as a composition of the
right-hand side of (\ref{eqspecialWSdef2}) with a monotone
function of $x$; this extension of (\ref{eqspecialWSdef2}) is
studied in \cite{CarrPCS}.
Therefore, in the present paper, we do \emph{not} focus on the
properties of stochastic processes $X$ which admit the particular weak
symmetry mapping (\ref{eqspecialWSdef2}). Instead, \emph{we show
in Theorem}~\ref{thmain} \textit{that any L\'evy process with one-sided jumps}
(\textit{subject to some regularity assumptions}) \textit{possesses a weak symmetry},
which is given by (\ref{eqMITsol}) and does not have to coincide
with (\ref{eqspecialWSdef2}).

\begin{remark}
Having established the weak reflection principle for a Markov process
$X$, with the weak symmetry mapping $\mathbf{W}^{\pm}$, one can
easily extend it to the process that arises as an independent
continuous stochastic time change of $X$, with the same $\mathbf
{W}^{\pm}$. This observation is made, for example, in \cite{CarrPCS}
for the particular choice of weak symmetry mapping (\ref
{eqspecialWSdef2}), but, of course, it remains valid for an
arbitrary $\mathbf{W}^{\pm}$.
\end{remark}

Next, let us demonstrate that the weak reflection principle provides
solutions to the problems outlined in Section~\ref{subseSRP},
for a class of processes that may not posses any strong symmetries.
Without loss of generality, we focus on the weak reflection principle
with an upper barrier. Namely, assuming that the process $X$, with
$X_0=0$, is strongly Markov, that it does not jump across the barrier
$U=0$ from below and that it possesses the upper weak symmetry $\mathbf
{W}^+$, we propose the following applications of our method.

First of all, the weak reflection principle allows us to solve the
static hedging problem in a model where $S=X^x$ is the underlying, with
$X^x_0=x\leq0$. Namely, for any admissible function $h$, with support
in $(-\infty,0)$, at any time up until and including the time when $S$
hits $0$, the price of a European-type option with maturity $T$ and the
payoff $h(S_T)-\mathbf{W}^+ h ( S_T  )$ is given by
\begin{eqnarray}\label{equopshformalpf}
&& \EE \bigl( h (S_T ) - \mathbf{W}^+ h (S_T ) \mid
\mathcal{F}_{t\wedge T_0} \bigr)\nonumber
\\
&&\qquad  = \EE \bigl( \EE \bigl( \bigl(h
(S_T ) - \mathbf{W}^+ h ( S_T ) \bigr)
\bone_{ \{T_0<T\}}\mid \mathcal{F}_{T\wedge T_0} \bigr) \mid \mathcal
{F}_{t\wedge T_0} \bigr) \nonumber
\\
&&\quad\qquad{}
+ \EE \bigl( \EE \bigl( \bigl(h (S_T ) - \mathbf{W}^+
h ( S_T ) \bigr) \bone_{ \{T_0\geq
T \}}\mid \mathcal{F}_{T\wedge T_0}
\bigr) \mid \mathcal{F}_{t\wedge T_0} \bigr)
\nonumber\\[-8pt]\\[-8pt]\nonumber
&&\qquad = \EE \bigl( \bigl(\EE h \bigl(X^s_{\tau} \bigr) - \EE
\mathbf{W}^+ h \bigl(X^s_{\tau} \bigr) \bigr)_{\tau=
T-T\wedge T_0, s=S_{T\wedge T_0}}
\bone_{ \{T_0<T \}} \mid \mathcal{F}_{t\wedge T_0} \bigr) \nonumber
\\
&&\qquad\quad{}
+ \EE \bigl( h (S_T ) \bone_{ \{T_0\geq
T \}} \mid
\mathcal{F}_{t\wedge T_0} \bigr)\nonumber
\\
&&\qquad
 = \EE \bigl( h (S_T )
\bone_{ \{\sup_{t\in
[0,T]}S_t<0 \}} \mid \mathcal{F}_{t\wedge T_0} \bigr),\nonumber
\end{eqnarray}
which coincides with the price of an up-and-out barrier option with the
terminal payoff function $h$. In the above, $X^s$ is a copy of the
original Markov process which is independent of $\mathcal{F}_{T\wedge
T_0}$. Then, the second equality follows from the strong Markov
property of $S$ and the fact that $\mathbf{W}^+ h$ is supported in
$[0,\infty)$. The last equality, in turn, follows from the continuity
of the running maximum of $S$ (which implies $S_{T_0}=0$) and the
definition of $\mathbf{W}^+ h$ as the image of $h$ under the weak
symmetry mapping.
Thus, in order to offset the risks associated with holding an
up-and-out barrier option with the terminal payoff function $h$ (i.e.,
hedge the barrier option), one needs to sell the European-type option
with the payoff function $h-\mathbf{W}^+ h$, and buy it back (at a
zero price) when, and if, the underlying hits~$0$. The static hedge
payoff corresponding to the up-and-out put option [i.e., with
$h(x)=(K-x)^+$] is computed in Section~\ref{seexample} for a
particular choice of L\'evy process $X$.

The weak reflection principle can also be used to express the joint
distribution of the process and its running maximum via the marginal
distribution of the process. Assuming again that $S=X^x$ with
$X^x_0=x\leq0$, we make use of (\ref{equopshformalpf}), to obtain,
for any $K<0$,
%
\begin{eqnarray}\label{eqjointlawcompdef}
\PP \Bigl( S_T\leq K, \sup_{t\in[0,T]}
S_t \geq0 \Bigr)
&=& \EE h_K(S_T) - \EE \bigl( h_K(S_T)
\bone_{ \{ \sup_{t\in
[0,T]} S_T < 0 \}} \bigr)
\nonumber\\[-8pt]\\[-8pt]\nonumber
&=&  \EE\mathbf{W}^+ h_K(S_T),
\end{eqnarray}
with $h_K = \bone_{(-\infty,K]}$.\footnote{There is a subtlety hidden in equation (\ref
{eqjointlawcompdef}). Namely, due to the discontinuity of $h_K$, its
weak symmetry image may be a \emph{generalized}, rather than
classical, function. Thus, the right-hand side of~(\ref
{eqjointlawcompdef}) should be understood as the action of $\mathbf
{W}^+ h_K$ on the density of $S_T$. Section~\ref{seexample} shows how
to make this interpretation rigorous.}
It is worth mentioning that, unlike the static hedging problem, the
computation of the joint law of a process does not require the
knowledge of the image of $h$ under the weak symmetry mapping. Indeed,
we only need to know the integral of $\mathbf{W}^+ h$ with respect to
the distribution of $S_T$. As a result, even though in certain cases it
may be advantageous to use (\ref{eqjointlawcompdef}) for the
computation of the joint law, this application does not fully utilize
the power of the weak reflection principle.
In Section~\ref{seexample}, we present the numerical implementation
of (\ref{eqjointlawcompdef}) for a particular choice of L\'evy
process $X$, and discuss the complexity of this method relative to the
existing ones; cf. \cite{s99,Kuznetsov3,Kuznetsov2,Kuznetsov1} and references therein.

Provided that the Markov process $X$ has a partial (integro-)
differential equation, or P(I)DE, associated with it (which is the case
for diffusions and L\'evy processes), the weak symmetry mapping
$\mathbf{W}^+$ allows us to solve the inverse problem associated with
this equation.
This problem is described briefly at the end of Section~\ref{subseSRP}, in the context of diffusion processes (Brownian motion,
in particular).
Note that the weak reflection principle allows us to solve this inverse
problem even when the associated equation does not possess any symmetries.
Section~\ref{seexample} provides a more detailed discussion of this
application for a particular choice of L\'evy process $X$.

\subsection{Prior results}


In view of the strong Markov assumption, it is natural to formulate the
problem of weak symmetry for jump-diffusions. In \cite{CarrN}, the
weak symmetry mapping is constructed for a class of diffusion processes
whose coefficients are only required to satisfy some regularity
conditions and do not have to be symmetric. Namely, consider a
diffusion process given by
\[
d X_t = \mu(X_t) \,dt + \sigma(X_t)
\,dB_t, %
\]
where $B$ is a Brownian motion. We assume that $\inf_{x\in\RR}\sigma
(x)>0$, the functions $\mu$ and $\sigma$ belong to $C^3(\RR)$, the
functions themselves and their first three derivatives have finite
limits at $-\infty$, and, for any $k=1,2,3$, the functions
$e^{(3-k)x}\llvert  \mu^{(k)}(x)\rrvert  $ and $e^{(3-k)x}\sigma^{(k)}(x)$ are
bounded over all $x>0$. Then Theorem~2.8 in \cite{CarrN} provides an
explicit integral representation for the weak symmetry mapping
associated with~$X$. For the sake of completeness, we present a
simplified corollary of this theorem here.

\begin{proposition}[(Carr and Nadtochiy \cite{CarrN})] Let $X$ be as above, and let
$h$ be a once weakly differentiable function, with support in $(-\infty
,0)$, such that its derivative is locally integrable and has a
modification with finite variation over $(-\infty,0)$. Then there
exists a continuous and exponentially bounded function $g$, with
support in $(0,\infty)$, such that
%
\begin{equation}
\label{eqgprop} \EE h(X_t) = \EE g(X_t)\qquad\mbox{for
all } t>0.
\end{equation}
Moreover, for any large enough $\gamma>0$, the function $g$ can be
computed as follows:
\begin{eqnarray*}
g(x) &=& \frac{2}{\pi i} \int_{\gamma- i\infty}^{\gamma+ i\infty}
\frac{w\psi_1(x,w)}{\partial_x\psi_1(0,w) - \partial_x\psi_2(0,w)}
\\
&&\hspace*{49pt}{}\times \int_{-\infty}^0 \frac{\psi_1(z,w)}{\sigma^2 (z )}
\exp \biggl(-2\int_0^{z} \frac{\mu(y)}{\sigma^2(y)} \,dy
\biggr) h(z) \,dz \,dw, %
\end{eqnarray*}
where $\psi_1$ and $\psi_2$ are the fundamental solutions of the
associated Sturm--Liouville equation
%
\begin{equation}
\label{eqSL} \frac{1}{2}\sigma^2(x) \frac{\partial^2}{\partial x^2}
\psi(x,w) + \mu(x)\frac{\partial}{\partial x} \psi(x,w) - w^2 \psi(x,w) = 0,
\end{equation}
determined uniquely, for all complex $w$ with large enough
$\operatorname{Re}(w)>0$, by the following conditions: $\psi_1(\cdot,w)$ is square integrable on $(-\infty,0)$, $\psi_2(\cdot,w)$ is
square integrable on $(0,\infty)$ and $\psi_1(0,w) = \psi_2(0,w) = 1$.
\end{proposition}

The above result shows that any regular enough diffusion possesses a
weak symmetry given by an explicit integral transform. Of course, in
order to implement this transform numerically, one needs to know the
fundamental solutions of the associated Sturm--Liouville equation.
These functions can be approximated efficiently by expanding them into
power series of $w$; cf. \cite{Titchmarsh}. Alternatively, one can
notice that if $\mu\equiv0$ and $\sigma$ is piecewise constant, then
$\psi_i$'s are piecewise linear-exponential (linear combinations of
exponentials). Thus we can approximate any function $\sigma$ with the
piecewise constant ones, then compute $\psi_i$'s in closed form, and
finally, obtain $g$ via numerical integration. Examples of functions
$g$ corresponding to a piecewise linear function $h$ can be found in
\cite{CarrN}.


The remainder of this paper is organized as follows. Section~\ref
{subseLaplace} introduces the notation, states the main assumptions
and formulates the weak symmetry problem for spectrally-negative L\'evy
processes. Section~\ref{subseest} contains the technical lemmas
needed for the proof of the main results, which are given in
Section~\ref{subsemain}. Theorem~\ref{thmain} provides the weak
symmetry mapping for spectrally-negative L\'evy processes, and
Corollary~\ref{cormainrate} addresses the computational aspects.
Finally, Section~\ref{seexample} illustrates the applications of the
weak reflection principle, for a particular choice of the L\'evy
process, and Section~\ref{sesummary} summarizes the results and
outlines the future research directions.

\section{Weak symmetry of spectrally negative L\'{e}vy processes}

\subsection{Problem formulation}\label{subseLaplace}
Consider a L\'evy process $(X_t)_{t\geq0}$, given by its initial
condition $X_0=0$ and the Laplace exponent $\psi$,
%
\begin{equation}
\label{eqLevyTriplet} \psi(\lambda) = \mu\lambda+ \frac{\sigma^2}{2}
\lambda^2 + \int_{-\infty}^0
\bigl(e^{\lambda x} - 1 - \lambda x \bigr) \Pi(dx),
\end{equation}
where $\Pi$ is the L\'evy measure of $X$, and
\[
\EE e^{\lambda X_t} = e^{t\psi(\lambda)}, %
\]
for all complex $\lambda$ for which both sides of the above equation
are well defined.
To make sure that the above expressions are well defined, at least, for
all $\lambda$ with positive real part, as well as to simplify some of
the derivations that follow, we make the following assumption on the L\'
evy triplet $(\mu,\sigma,\Pi)$.

\begin{ass}\label{ass1}
We assume that $\mu\in\RR$, $\sigma>0$ and that $\Pi$ is a $\sigma
$-finite Borel measure on $(-\infty,0)$, satisfying
\[
\int_{-\infty}^{0} x^2 \Pi(dx) < \infty
\quad\mbox{and}\quad\int_{-\infty}^{-1} \llvert x\rrvert
e^{-\zeta x} \Pi(dx) < \infty, %
\]
with some $\zeta\geq0$.
\end{ass}

The process $X$ is called \emph{spectrally negative} because it is
only allowed to have negative (i.e., downward) jumps; cf. \cite
{Kuznetsov1}. The reason for such a restriction is that the process
must not jump across the upper barrier in order for the weak reflection
principle to hold. Of course, in the case of a lower barrier, one needs
to consider $X$ with positive jumps only.

Our goal is to construct a weak symmetry mapping for the process $X$.
Namely, for any given admissible function $h\dvtx  \RR\rightarrow\RR$,
with $\operatorname{supp}(h)\subset(-\infty,0)$, we would like to
find a measurable function $g\dvtx  \RR\rightarrow\RR$, with
$\operatorname{supp}(g)\subset[0,\infty)$, such that
%
\begin{equation}
\label{eqweaksymprobdef1} \EE h(X_t) = \EE g(X_t)\qquad\forall t>0.
\end{equation}
To ensure that the expectation of $h(X_t)$ is well defined, we need to
make some additional assumptions on $h$ and $\Pi$.

\begin{ass}\label{ass2}
We assume that $\operatorname{supp}(h)\subset(-\infty,0)$ and that
there exists $\hat{h}\in\mathbb{L}^1(\RR)$, such that the function
$x\mapsto e^{\zeta x}h(x)$, defined for all $x\in\RR$, is a Fourier
transform of $\hat{h}$ (with $\zeta$ given in Assumption~\ref{ass1}).\footnote{We use the convention that Assumptions~\ref{ass1}
and~\ref{ass2} hold for all the subsequent derivations, except the
statements of lemmas, theorems and corollaries, where all assumptions
are stated explicitly.}
\end{ass}

Taking the Laplace transform with respect to $t$ on both sides of (\ref
{eqweaksymprobdef1}), we obtain
%
\begin{equation}
\label{eqinvpr1} \EE\int_0^{\infty} e^{-\lambda t}
h(X_t) \,dt = \EE\int_0^{\infty}
e^{-\lambda t} g(X_t) \,dt.
\end{equation}
Assumption~\ref{ass2} implies that
%
\begin{eqnarray}\label{eqEhLaplaceest}
\EE\int_0^{\infty} e^{-\lambda t}
\bigl\llvert h(X_t)\bigr\rrvert \,dt &\leq& \int_0^{\infty}
e^{-\lambda t} \EE e^{-\zeta X_t} \,dt
\nonumber\\[-8pt]\\[-8pt]\nonumber
&=& \int_0^{\infty}e^{-(\lambda-\psi(-\zeta)) t} \,dt <\infty
\end{eqnarray}
holds for all real $\lambda>\psi(-\zeta)\vee0$.
Due to the uniqueness of the Laplace inverse, our problem is equivalent
to the following: find a measurable function $g\dvtx  \RR\rightarrow\RR$,
with $\operatorname{supp}(g)\subset[0,\infty)$, such that (\ref
{eqinvpr1}) holds for all large enough $\lambda>0$.

Next, recall that the Laplace transform in time of the expectation of a
function of a Markov process is given by the value of an integral
operator, called the \emph{resolvent operator}, applied to this
function. In particular, as follows, for example, from Theorem 2.7 in
\cite{Kuznetsov1}, for the spectrally negative L\'evy process $X$, we have
\[
\EE\int_0^{\infty} e^{-\lambda t}
f(X_t) \,dt = \int_{\RR} \bigl(
\Phi'(\lambda)e^{-\Phi(\lambda)x}-W^{\lambda
}(-x) \bigr) f(x) \,dx,
\]
which holds for any measurable function $f$ for which the integral on
the left-hand side is absolutely convergent, with
\[
\Phi(\lambda):=\sup\bigl\{q \geq0\dvtx  \psi(q)=\lambda\bigr\},\qquad \lambda\geq0
\]
and with $W^{\lambda}$ being the \emph{$\lambda$-scale function} of
$X$. Recall that, for every real $\lambda\geq0$, the $\lambda$-scale
function $W^{\lambda}\dvtx \RR\rightarrow[0,\infty)$ is uniquely defined
as the right-continuous function that takes value zero in $(-\infty,0)$ and satisfies
\[
\int_{\RR} e^{-w x} W^{\lambda}(x)\,dx =
\frac{1}{\psi(w)-\lambda}, %
\]
for all $w\geq\Phi(\lambda)$. For further details on the theory of
scale functions, we refer the reader to \cite{Kuznetsov1} and the
references therein.
Thus (\ref{eqinvpr1}) is equivalent to
\[
\int_0^{\infty} \Phi'(
\lambda)e^{-\Phi(\lambda)x} g(x) \,dx = \Upsilon(\lambda), %
\]
where
%
\begin{equation}
\label{eqUpsdef} \qquad\Upsilon(\lambda):= \EE\int_0^{\infty}
e^{-\lambda t} h(X_t) \,dt = \int_0^{\infty}
\bigl( \Phi'(\lambda)e^{\Phi(\lambda
)x}-W^{\lambda}(x) \bigr)
h(-x) \,dx.
\end{equation}
Notice that $\psi$ is continuous and strictly increasing on $
(\llvert  \mu\rrvert  /\sigma^2,\infty )$, exploding at infinity.
Hence, we can change the variables in the above equations to obtain an
equivalent formulation of the problem. Namely, we search for a
measurable function $g$, such that, for all large enough real $\lambda
>0$, the following holds:
%
\begin{equation}
\label{eqmirrorImagedef} \int_0^{\infty} e^{-\lambda x} g(x)
\,dx = \psi'(\lambda) \Upsilon\bigl(\psi(\lambda)\bigr),
\end{equation}
with
%
\begin{equation}
\label{eqrhsdef} \psi'(\lambda) \Upsilon\bigl(\psi(\lambda)\bigr) =
\int_0^{\infty} \bigl( e^{\lambda x} -
\psi'(\lambda) W^{\psi
(\lambda)}(x) \bigr) h(-x) \,dx.
\end{equation}
%

Notice that the weak symmetry problem (\ref{eqmirrorImagedef}) now
looks exactly like the Laplace transform inversion. However, there is a
major difference. In the classical problem of inverting a Laplace
transform, we typically know that the right-hand side is a Laplace
transform of some function, and we need to find a mapping that would
recover this function. In the present case, we need to prove that $\psi
'(\lambda) \Upsilon(\psi(\lambda))$ is, indeed, a Laplace transform
of some function $g$, which is not obvious a priori. In addition, we
need to propose a method to recover $g$. Of course, there exist several
methods for inverting the Laplace transform; cf. \cite{WidderLaplace,Davies}. However, the conditions that are required for some of
these methods to succeed are expressed through the original function
(in our case, $g$), rather than the transformed one [in the present
case, $\psi'(\lambda) \Upsilon(\psi(\lambda))$]. Since, a priori,
we know very little about function $g$ (e.g., we do not even know if it
exists), we cannot apply any of the existing results on the Laplace
transform inversion to solve the weak symmetry problem. Instead, we
will show, by hand, that the desired function $g$ exists and that it
can be recovered via the classical Bromwich integral; cf. \cite{Davies}.

\subsection{A priori estimates}\label{subseest}
In this subsection, we establish the analytic continuation of $\psi
'(\lambda) \Upsilon(\psi(\lambda))$ to a complex half plane of the
form $H_{R}= \{ w\dvtx  \Re(w)>R  \}$ and provide some useful
estimates of its absolute value.

\begin{remark}\label{remrem1}
Notice that $\Upsilon(\lambda)$ can be easily extended to a half
plane $H_R$, via its probabilistic representation given by the first
identity in (\ref{eqUpsdef}). However, under the change of
variables $\lambda\mapsto\psi(\lambda)$, the half plane transforms
into a smaller domain which is not sufficient for our purposes. To
obtain an analytic extension of $\Upsilon(\psi(\lambda))$, one would
need the probabilistic representation in (\ref{eqUpsdef}) to hold
in a domain where the real part of $\lambda$ is unbounded from below,
which is usually impossible. More precisely, the conditions that
function $h$ has to satisfy, in order for the probabilistic
representation in (\ref{eqUpsdef}) to hold in the desired domain,
are extremely restrictive and rather implicit. For example, if $X$ is a
martingale, these conditions exclude all convex functions $h$, except
zero. In addition, even for those functions $h$ for which the
probabilistic representation in (\ref{eqUpsdef}) holds in the
desired domain (although we do not know how to characterize this set
explicitly), the standard estimates do not provide sufficient
information about the asymptotic behavior of $\Upsilon(\psi(\lambda
))$, as $\llvert \lambda\rrvert  \rightarrow\infty$, which is needed to solve the
weak symmetry\vspace*{1pt} problem (\ref{eqmirrorImagedef}).

Alternatively, one might be tempted to use the analytic continuation of
$W^{\psi(\lambda)}$ to extend $\Upsilon(\psi(\lambda))$ to a
complex half plane via (\ref{eqrhsdef}). Indeed, it is well known
(cf.~\cite{Kuznetsov1}) that $W^{\lambda}$ can be extended
analytically to the entire complex plane. However, to the best of our
knowledge, there exist no estimates of this extension [more precisely,
we need an estimate of $e^{\lambda x} - \psi'(\lambda) W^{\psi
(\lambda)}(x)$], which, in particular, would guarantee that the
integral on the right-hand side of (\ref{eqrhsdef}) is well
defined for all $\lambda$ in a half plane $H_{R}$. In fact, it is easy
to construct an example of a L\'evy process whose infinitesimal
generator has a nontrivial spectrum, and hence the integral on the
right-hand side of (\ref{eqrhsdef}) is not well defined for some
$\lambda$.

The difficulties described above explain why we are forced to construct
the analytic continuation of $\Upsilon(\psi(\lambda))$, and
investigate its asymptotic behavior, by hand.
\end{remark}

It turns out that Fourier transform offers a natural way to obtain the
desired analytic continuation.
Recall that, due to Assumption~\ref{ass2},
\[
h(x) = e^{-\zeta x} \int_{\RR} e^{-ixz}
\hat{h}(z) \,dz. %
\]
%
Making use of (\ref{eqEhLaplaceest}), we
apply Fubini's theorem to obtain
%
\begin{eqnarray}
\label{eqrep1}
&& \psi'(\lambda) \Upsilon\bigl(\psi(\lambda)\bigr)\nonumber
\\
&&\qquad =
\int_0^{\infty} \bigl( e^{\lambda x} -
\psi'(\lambda) W^{\psi
(\lambda)}(x) \bigr) h(-x) \,dx
\\
&&\qquad = \int_{\RR} \int_0^{\infty}
\bigl( e^{\lambda x} - \psi'(\lambda ) W^{\psi(\lambda)}(x) \bigr)
e^{(\zeta+iz)x} \,dx\, \hat{h}(z) \,dz, \nonumber
\end{eqnarray}
for all large enough real $\lambda>0$.
Our next goal is to extend the above representation to a complex half
plane $H_R$ and estimate its absolute value from above.
Let us analyze the inner integral on the right-hand side of (\ref
{eqrep1}). Notice that it can be viewed as
\[
\int_0^{\infty} \bigl( e^{\lambda x} -
\psi'(\lambda) W^{\psi
(\lambda)}(x) \bigr) e^{-wx} \,dx,
\]
evaluated at $w= -\zeta- iz$.
For $\lambda>\llvert  \mu\rrvert  /\sigma^2$ and $w\in(\lambda,\infty)$, the
above integral can be computed explicitly. Using the definition of the
scale function $W^{\psi(\lambda)}$ (cf. \cite{Kuznetsov1}), we obtain
%
\begin{eqnarray}\label{eqmirror2}
&& \int_0^{\infty} \bigl(
e^{\lambda x} - \psi'(\lambda) W^{\psi
(\lambda)}(x) \bigr)
e^{-wx} \,dx\nonumber
\\
&&\qquad  = \int_0^{\infty}
e^{\lambda x - wx} \,dx - \psi'(\lambda) \int_0^{\infty}
e^{-wx} W^{\psi(\lambda)}(x) \,dx
\nonumber\\[-8pt]\\[-8pt]\nonumber
&&\qquad = \frac{1}{w-\lambda} - \frac{\psi'(\lambda)}{\psi(w) - \psi
(\lambda)}
\\
&&\qquad = \frac{\psi(w) - \psi(\lambda) - \psi'(\lambda)(w-\lambda
)}{(w-\lambda)^2}
\frac{w-\lambda}{\psi(w) - \psi(\lambda)}. \nonumber
\end{eqnarray}
The right-hand side of the above is analytic everywhere in $w\in
H_{-\zeta}= \{ w\dvtx  \Re(w)>-\zeta \}$, except, possibly,
the zeros of $\psi(w) - \psi(\lambda)$, where it has poles (on every
compact, there is at most a finite number of such points).
The left-hand side of the above is analytic in $H_{-\zeta}$ (we can
differentiate with respect to $w$ inside the integral, and the
resulting expression is bounded, uniformly over $w$, by an absolutely
integrable function of $x$). Therefore, the right-hand side, in fact,
does not have any poles in $H_{-\zeta}$, and the above equality holds
for all $w\in H_{-\zeta}$. By continuity, the equality holds at
$w=-\zeta-iz$
%
\begin{eqnarray}\label{eqLevyfinal1}
&& \int_0^{\infty} \bigl(
e^{\lambda x} - \psi'(\lambda) W^{\psi
(\lambda)}(x) \bigr)
e^{(\zeta+iz)x} \,dx
\nonumber\\[-8pt]\\[-8pt]\nonumber
&&\qquad = \frac{\psi'(\lambda)}{\psi(\lambda) - \psi(-\zeta-iz)} - \frac
{1}{\lambda+ \zeta+ iz},
\end{eqnarray}
for all large enough real $\lambda>0$.
Next, we need to estimate the absolute value of the right-hand side of
the above.
To do this, we will use the following lemmas, which are also crucial
for the proofs of the main result. These lemmas describe the asymptotic
behavior of $1/(\psi(\lambda) - \psi(-\zeta-iz))$, and they can be
viewed as the core technical result of the paper. Their proofs are
given in Appendix~\ref{appA}.


\begin{lemma}\label{leIm1}
Let Assumption~\ref{ass1} hold. For all $u,z\in\RR$ and $v>0$, we have
\begin{eqnarray*}
\Im \bigl(\psi(v + iu) - \psi(-\zeta-iz) \bigr)
&=& \mu(u+z) + (v u - z \zeta)
\sigma^2 + \overline{\overline{o}}(uv)
\\
&&{} + \operatorname{sign}(u) \bigl
\llvert \overline{\overline{o}}\bigl(u^2\bigr) \bigr\rrvert +
\operatorname{sign}(z) \bigl\llvert \overline{\overline{o}}\bigl(z^2
\bigr) \bigr\rrvert + \overline{\overline{o}}(z), %
\end{eqnarray*}
where we denote by $\overline{\overline{o}} (f(u,v,z) )$
any function of $(u,v,z)$, such that $\overline{\overline{o}}
(f(u,v,\break z) )/ f(u,v,z)$ is absolutely bounded over all $u,z\in\RR
$ and $v>0$, and it converges to zero, as $\llvert  f(u,v,z)\rrvert  \rightarrow\infty$.
\end{lemma}

\begin{lemma}\label{leIm2}
Let Assumption~\ref{ass1} hold. For all $u,z\in\RR$ and $v>0$, we have
\begin{eqnarray*}
\Im \bigl(\psi(v + iu) - \psi(-\zeta-iz) \bigr) &=& \mu(u+z) + (v u - z \zeta)
\sigma^2 + \overline{\overline{o}}(uv) + \overline{\overline{o}}
\bigl(z^2-u^2\bigr)
\\
&&{}  + \operatorname{sign}(z+u) \bigl
\llvert \overline{\overline{o}} \bigl((z+u)^2 \bigr)\bigr\rrvert +
\overline{\overline{o}}(z), %
\end{eqnarray*}
where $\overline{\overline{o}} (f(u,v,z) )$ has the same
meaning as in Lemma~\ref{leIm1}.
\end{lemma}

\begin{lemma}\label{leRe}
Let Assumption~\ref{ass1} hold. For all $u,z\in\RR$ and $v>0$, we have
\begin{eqnarray*}
\Re \bigl(\psi(v + iu) - \psi(-\zeta-iz) \bigr) &=& \mu(v+\zeta) +
\bigl(v^2 + z^2 - u^2 \bigr)
\sigma^2/2 + \overline{\overline{o}}\bigl(z^2-u^2
\bigr)
\\
&&{} + \bigl\llvert \overline{\overline{o}}(uv)\bigr\rrvert + \bigl\llvert
\overline {\overline{o}}\bigl(v^2\bigr)\bigr\rrvert + \overline{
\overline{o}}(z) + \overline{\overline{o}}(1), %
\end{eqnarray*}
where $\overline{\overline{o}} (f(u,v,z) )$ has the same
meaning as in Lemma~\ref{leIm1}.
\end{lemma}

\begin{lemma}\label{lepsiest1}
Let Assumption~\ref{ass1} hold. For any $\varepsilon\in(0,1)$,
there exist $R_1>0$, $R_2>0$, $c_1>0$ and $c_2>0$, such that the inequality
\begin{eqnarray*}
\bigl\llvert \psi(v + iu) - \psi(-\zeta-iz)\bigr\rrvert &\geq& c_1 v
\bigl(\llvert z\rrvert + v\bigr) \bone_{ \{\llvert  \llvert  u\rrvert   - \sqrt{z^2 + v^2} \rrvert
\leq\varepsilon v \}}
\\
&&{} + c_2 \bigl
\llvert u^2 - \bigl(z^2 + v^2\bigr) \bigr
\rrvert \bone_{ \{\llvert  \llvert  u\rrvert   -
\sqrt{z^2 + v^2} \rrvert   > \varepsilon v \}} %
\end{eqnarray*}
holds for all $u,v,z\in\RR$, with $v>R_1$ and $\llvert  z\rrvert  >R_2$.
\end{lemma}

\begin{lemma}\label{lepsiest2}
Let Assumption~\ref{ass1} hold. There exists a constant $R_1>0$, such
that, for any $v>R_1$, there exist $R_2=R_2(v)>0$, $R_3=R_3(v)>0$ and
$c=c(v)>0$, such that the inequality
\[
\bigl\llvert \psi(v + iu) - \psi(-\zeta-iz)\bigr\rrvert ^2 \geq c
\bigl( \bigl(z^2 - u^2 \bigr)^2 +
z^2 \bigr), %
\]
holds for all $u,z\in\RR$, with $\llvert  z\rrvert  >R_2$ and $\llvert  u\rrvert  >R_3$.
\end{lemma}

The following lemma makes use of the above results to show that the
right-hand side of (\ref{eqLevyfinal1}) is absolutely integrable
with respect to $\hat{h}(z)\,dz$, even for complex $\lambda$.

\begin{lemma}\label{lepsiprimeUpsilonest}
Let Assumption~\ref{ass1} hold. Then there exist constants $c$ and
$R_1$, such that the following holds, for all $u\in\RR$, all $v\geq
R_1$ and all $\hat{h}\in\mathbb{L}^1\cup\mathbb{L}^2$: $\psi(v +
iu) \neq\psi(-\zeta-iz)$, for all $z\in\RR$ and
%
\begin{eqnarray}\label{eqcor2est1}
&& \int_{\RR} \biggl\llvert \frac{\psi'(v + iu)}{\psi(v + iu) - \psi
(-\zeta-iz)} -
\frac{1}{v + iu + \zeta+ iz} \biggr\rrvert \bigl\llvert \hat{h}(z)\bigr\rrvert \,dz
\nonumber\\[-8pt]\\[-8pt]\nonumber
&&\qquad \leq c \llvert v + iu\rrvert \bigl(\llVert \hat{h}\rrVert _{\mathbb{L}^1(\RR)}\wedge\llVert
\hat {h}\rrVert _{\mathbb{L}^2(\RR)} \bigr).
\end{eqnarray}
\end{lemma}

\begin{pf}
Fix arbitrary $\varepsilon\in(0,1)$. Lemma~\ref{lepsiest1} yields
that there exist $R_1,R_2>0$ and $c_i$'s, such that the following
estimates hold, for all $u\in\RR$, all $v\geq R_1$ and all $\llvert  z\rrvert  \geq R_2$:
%
\begin{eqnarray}\label{eqLevyfinaloneoverpsi}
\qquad&& \biggl\llvert \frac{1}{\psi(v+iu) - \psi(-\zeta-iz)}\biggr\rrvert ^2 \nonumber
\\
&&\qquad  \leq \frac{c_1}{v^2  (\llvert  z\rrvert   + v )^2} \bone_{ \{\llvert   \llvert  u\rrvert   - \sqrt{z^2 + v^2} \rrvert   \leq
\varepsilon v \}} + \frac{c_2}{\varepsilon^2 v^2 (v^2 + z^2)}
\bone_{ \{\llvert   \llvert  u\rrvert   - \sqrt{z^2 + v^2} \rrvert   > \varepsilon
v \}}
\\
&&\qquad  \leq \frac{c_3}{(1+\llvert  z\rrvert  )^2}.\nonumber
\end{eqnarray}
It is also easy to see [by a direct examination of (\ref
{eqLevyTriplet})] that, for all large enough $\llvert \lambda\rrvert  $, such that
$\Re(\lambda)\geq-\zeta$, we have
%
\begin{equation}
\label{eqpsipsiprimeasymp} \bigl\llvert \psi''(\lambda)\bigr
\rrvert \leq c_4,\qquad
\bigl\llvert \psi'(\lambda)\bigr
\rrvert \leq c_{5} \llvert \lambda\rrvert,\qquad c_{6} \llvert
\lambda\rrvert ^2\leq\bigl\llvert \psi(\lambda)\bigr\rrvert \leq
c_{7} \llvert \lambda\rrvert ^2,
\end{equation}
with some strictly positive constants $c_i$.
Collecting the above and possibly increasing $R_1$, we obtain
\begin{eqnarray*}
&& \int_{\RR} \biggl\llvert \frac{\psi'(v + iu)}{\psi(v + iu) - \psi
(-\zeta-iz)} - \frac{1}{v + iu + \zeta+ iz} \biggr\rrvert \bigl\llvert \hat{h}(z)\bigr\rrvert \,dz
\\
&&\qquad \leq \frac{c_8\llvert  v + iu\rrvert  }{\llvert  v + iu\rrvert  ^2 - R_2^2} \int_{\llvert  z\rrvert  \leq R_2} \bigl\llvert \hat {h}(z)
\bigr\rrvert \,dz
\\
&&\quad\qquad{}+ \int_{\llvert  z\rrvert  >R_2} \biggl\llvert \frac{\psi'(v + iu)}{\psi(v + iu) - \psi
(-\zeta-iz)}
\biggr\rrvert \bigl\llvert \hat{h}(z)\bigr\rrvert \,dz %
\\
&&\quad\qquad{}
+ \int_{\RR} \biggl\llvert \frac{1}{v + iu + \zeta+ iz} \biggr\rrvert
\bigl\llvert \hat {h}(z)\bigr\rrvert \,dz %
\\
&&\qquad \leq \frac{c_8\llvert  v + iu\rrvert  }{\llvert  v + iu\rrvert  ^2 - R_2^2} \int_{\llvert  z\rrvert  \leq R_2} \bigl\llvert \hat {h}(z)
\bigr\rrvert \,dz + c_9 \llvert v + iu\rrvert \int_{\llvert  z\rrvert  >R_2}
\frac{1}{1+\llvert  z\rrvert  } \bigl\llvert \hat{h}(z)\bigr\rrvert \,dz %
\\
&&\quad\qquad{} + c_{10}\int_{\RR}\frac{1}{v + \zeta+ \llvert  z+u\rrvert  }\bigl\llvert
\hat{h}(z)\bigr\rrvert \,dz, %
\end{eqnarray*}
which yields (\ref{eqcor2est1}), after an application of the Cauchy
inequality.
\end{pf}

Recall that, due to equations (\ref{eqrep1}) and (\ref
{eqLevyfinal1}), the representation
%
\begin{equation}
\label{eqanExt} \psi'(\lambda)\Upsilon\bigl(\psi(\lambda)\bigr) =
\int_{\RR} \biggl( \frac{\psi'(\lambda)}{\psi(\lambda) - \psi
(-\zeta-iz)} - \frac{1}{\lambda+ \zeta+ iz}
\biggr) \hat{h}(z) \,dz
\end{equation}
is well defined and holds for all large enough real $\lambda>0$.
Lemma~\ref{lepsiprimeUpsilonest} shows that the right-hand side of
the above is well defined for all $\lambda\in H_{R_1}= \{\lambda\dvtx  \Re(\lambda)>R_1 \}$, with $R_1>0$ given in Lemma~\ref
{lepsiprimeUpsilonest}. In fact, it is easy to deduce that the
right-hand side of (\ref{eqanExt}) is analytic in $H_{R_1}$.
To see this, first notice that the integrand in (\ref{eqanExt}) is
analytic in $\lambda\in H_{R_1}$ (as the denominators cannot vanish).
Then differentiate, formally, inside the integral, and apply the same
estimates as in the proof of Lemma~\ref{lepsiprimeUpsilonest} to
show that the integral of the derivative is absolutely convergent, for
any $\lambda\in H_{R_1}$.
Thus we have proved the following corollary.

\begin{cor}\label{coranExt}
Let Assumptions~\ref{ass1} and~\ref{ass2} hold, and let $R_1>0$ be
the constant appearing in Lemma~\ref{lepsiprimeUpsilonest}. Then the
function $\lambda\mapsto\psi'(\lambda)\Upsilon(\psi(\lambda))$
[defined in (\ref{eqrhsdef}), for all large enough real $\lambda
>0$] can be extended analytically to $H_{R_1}= \{\lambda\dvtx
\Re(\lambda)>R_1 \}$ via (\ref{eqanExt}).
\end{cor}

\subsection{Main results}\label{subsemain}
Now, we have everything we need to solve the weak symmetry problem
(\ref{eqmirrorImagedef}).
First, for any $r>0$, we introduce the following function of $x\in\RR$:
%
\begin{eqnarray}\label{eqgrdef}
g_r(x) &=& \frac{1}{2\pi} \int_{-r}^r
e^{(\gamma+iu) x}
\int_{\RR} \biggl( \frac{\psi'(\gamma+iu)}{\psi(\gamma+iu) -
\psi(-\zeta-iz)}
\nonumber\\[-8pt]\\[-8pt]\nonumber
&&\hspace*{121pt}{} -\frac{1}{\gamma+iu + \zeta+ iz} \biggr) \hat{h}(z) \,dz \,du,
\end{eqnarray}
%
with a large (but fixed) constant $\gamma>0$.
To ensure that $g_r$ is well defined, we assume that $\gamma>R_1$,
where $R_1$ is the constant appearing in Lemma~\ref{lepsiprimeUpsilonest}.
In this section, we show that $g_r$ has a limit $g$, as $r\rightarrow
\infty$, and that $g$ is the weak symmetry image of~$h$ (as discussed
in Section~\ref{subseWRP}).
However, before we present the main result of the paper, we need to
relax our assumptions on function $h$.
Notice that Assumption~\ref{ass2} excludes some functions $h$ that
are important for applications, such as the indicator
functions.\footnote{We thank the anonymous referee for pointing this
out.} Indeed, if $h=\bone_{(-\infty,K]}$, then
\[
\hat{h}(z) = \frac{\exp(K(\zeta+ iz))}{2\pi(\zeta+ iz)}, %
\]
which is not absolutely integrable over $\RR$. Nevertheless, the above
function belongs to $\mathbb{L}^2(\RR)$, and the right-hand side of
(\ref{eqgrdef}) is well defined for all $\hat{h}\in\mathbb
{L}^1(\RR)\cup\mathbb{L}^2(\RR)$, due to Lemma~\ref
{lepsiprimeUpsilonest}. Thus we extend the scope of our analysis to
include all functions $h$ that satisfy the following assumption (which
is a strictly weaker version of Assumption~\ref{ass2}).

\begin{ass}\label{ass3}
We assume that $\operatorname{supp}(h)\subset(-\infty,0)$ and that
there exists $\hat{h}\in\mathbb{L}^1(\RR)\cup\mathbb{L}^2(\RR)$,
such that the function $x\mapsto e^{\zeta x}h(x)$, defined for all
$x\in\RR$, coincides almost everywhere with the Fourier transform of
$\hat{h}$ (with $\zeta$ given in Assumption~\ref{ass1}).
\end{ass}

Recall that, in order to construct the weak symmetry image of $h$, we
need to consider the expectation of $h(X_t)$. However, if $\hat{h}\in
\mathbb{L}^2(\RR)$, we cannot guarantee that the expectation of
$h(X_t)$ is well defined: in this case, $h$ may not be locally bounded.
Nevertheless, if $X_t$ has a well-behaved density $p_t$ [such that
$e^{-\zeta x}p_t(x)$ is square integrable over $x\in\RR$], the
expectation of $h(X_t)$ is well defined since the associated integral
\[
\EE h(X_t) = \int_{\RR} \bigl(e^{\zeta x}h(x)
\bigr) \bigl(e^{-\zeta x}p_t(x) \bigr) \,dx %
\]
is absolutely convergent, which follows from the Cauchy inequality.
The following discussion shows that $p_t$ does possess the desired
properties and in addition, provides some auxiliary constructions
needed to formulate the main result.

Denote by $\mathcal{D}$ the space of all functions $f\dvtx \RR\rightarrow
\RR$, which are $3$ times continuously differentiable and satisfy
\[
\sup_{x\in\RR} \bigl\llvert e^{kx} f^{(n)}(x)
\bigr\rrvert < \infty,\qquad n=0,1,2,3, k=0,1,2\ldots.
\]
We equip $\mathcal{D}$ with the topology generated by the above family
of semi-norms and consider $\mathcal{D}^*$---the dual of $\mathcal
{D}$, consisting of all continuous linear functionals on~$\mathcal{D}$.
Using Fubini's theorem, it is easy to see that $g_r \in\mathcal
{D}^*$, for any $r>0$. In addition, Assumption~\ref{ass1} implies
that the marginal density $p_t$ of $X_t$ is well defined, for all \mbox{$t>0$}.
For any fixed $t>0$ and $k=0,1,2,\ldots,$ the Fourier transform of
$x\mapsto e^{kx}p_t(x)$ is given by
%
\begin{equation}
\label{eqfourdens} \int_{\RR} e^{-iux}e^{kx}p_t(x)
\,dx = \EE e^{(k-iu)X_t} = \exp \bigl(t\psi(k-iu) \bigr),
\end{equation}
which can be computed for all $u\in\RR$ as an analytic continuation
of its values in the domain $u\in \{-iy \mid y\in(k,\infty
) \}$, where the integral on the left-hand side is a priori known
to be absolutely convergent. Applying standard estimates to the
integral term in $\psi(k+iu)$ [cf. (\ref{eqLevyTriplet})], we can
easily deduce that
%
\begin{equation}
\label{eqfourdensasymp} \psi(k+iu) \sim-c u^2,\qquad \llvert u\rrvert \rightarrow
\infty,
\end{equation}
where $c$ is a positive constant. This, in particular, shows that the
Fourier transform of $x\mapsto e^{kx}p_t(x)$ belongs to the Schwartz
space (i.e., the space of infinitely smooth functions, decaying at $\pm
\infty$ faster than any power). 
Using the standard properties of Fourier transform, we conclude that
$x\mapsto e^{kx}p_t(x)$ belongs to the Schwartz space as well.
In addition, we notice that relations (\ref{eqfourdens})--(\ref{eqfourdensasymp}) hold for all $k\geq-\zeta$,
due to Assumption~\ref{ass1}.
Thus, we have proved the following lemma.

\begin{lemma}\label{leptl2}
Let\vspace*{1pt} Assumption~\ref{ass1} hold. For any $t>0$, we have $p_t\in
\mathcal{D}$, and the mapping $x\mapsto e^{-\zeta x}p_t(x)$ belongs to
$\mathbb{L}^1(\RR)\cap\mathbb{L}^2(\RR)$.
\end{lemma}

Finally, we present the main result of this paper.

\begin{theorem}\label{thmain}
Let Assumptions~\ref{ass1} and~\ref{ass3} hold, and let $g_r$
be given by (\ref{eqgrdef}), with an arbitrary $\gamma>R_1$, where
$R_1$ is the constant appearing in Lemma~\ref{lepsiprimeUpsilonest}.
Then, the following holds.
\begin{itemize}
\item As $r\rightarrow\infty$, $g_r$ converges weakly to a
generalized function $g\in\mathcal{D}^*$, which has support in
$[0,\infty)$ and satisfies
\[
\label{eqmirrorImdefweak} \langle g, p_t \rangle= \EE h(X_t)
\qquad\mbox{for all } t>0.
\]
Moreover, there exist constants $c_1, c_2>0$, independent of $\hat
{h}$, $r$ and $t$, such that, for all $r>0$ and all $t>0$, we have
%
\begin{equation}
\label{eqmainrategrweakpt} \bigl\llvert \langle g- g_{r}, p_t\rangle
\bigr\rrvert \leq c_1 \bigl(\llVert \hat{h}\rrVert _{\mathbb{L}^1(\RR)}
\wedge\llVert \hat{h}\rrVert _{\mathbb{L}^2(\RR)} \bigr) \int_{r}^{\infty}
u \exp \bigl(-c_2 t u^2 \bigr) \,du.
\end{equation}

\item If, in addition, Assumption~\ref{ass2} holds, then the
restriction of $g$ to the interval $(0,\infty)$ coincides with a
continuous function, which has at most exponential growth at infinity
[i.e., $\llvert  g(x)\rrvert  $ is bounded by a constant times an exponential, for all\vspace*{1pt}
large enough values of $x$]. Moreover, there exists a constant $c>0$,
independent of $\hat{h}$, $r$ and $x$, such that, for all large enough
$r>0$ and all $x>0$, we have
%
\begin{equation}
\label{eqmainrategr} \bigl\llvert g(x) - g_r(x)\bigr\rrvert \leq c
\frac{e^{\gamma x}}{x} \biggl(\frac{1}{r}\llVert \hat{h}\rrVert
_{\mathbb{L}^1(\RR)} + \int_{\llvert  z\rrvert  >r/2} \bigl\llvert \hat{h}(z)\bigr
\rrvert \,dz \biggr).
\end{equation}

\item If, in addition, Assumption~\ref{ass2} holds and
%
\begin{equation}
\label{eqtailcond} \int_0^{\infty} \biggl(\int
_{\llvert  z\rrvert  >\llvert  u\rrvert  } \bigl\llvert \hat{h}(z)\bigr\rrvert \,dz
\biggr)^2 \,du<\infty,
\end{equation}
then $g$ is locally integrable in $\RR$ and continuous in $\RR
\setminus \{0  \}$, with at most exponential growth at
infinity, and
%
\begin{equation}
\label{eqmirrorImdefstrong} \EE g(X_t) = \EE h(X_t)\qquad\mbox{for
all } t>0.
\end{equation}
\end{itemize}
\end{theorem}

\begin{remark}
Notice that the uniqueness of function $g$, having support in
$[0,\infty)$ and satisfying (\ref{eqmirrorImdefstrong}), follows
from the uniqueness of the Laplace inverse and from equation (\ref
{eqmirrorImagedef}), derived in Section~\ref{subseLaplace}.
\end{remark}

The proof of the above theorem is given in Appendix~\ref{appB}.
Theorem~\ref{thmain} shows that any spectrally negative L\'evy
process $X$, satisfying Assumption~\ref{ass1}, possesses an upper
weak symmetry, with the space of test functions $\mathcal{B}^-_0$
consisting of all functions that satisfy Assumption~\ref{ass2}. The
associated weak symmetry image transformation is given by
%
\begin{eqnarray}\label{eqMITsol}
\mathbf{W}^+ h(x) &=& g(x)\nonumber
\\
&=& \frac{1}{2\pi} \int_{\RR} e^{(\gamma
+iu) x}
\int_{\RR} \biggl(
\frac{\psi'(\gamma+iu)}{\psi(\gamma+iu) -
\psi(-\zeta-iz)}
\\
&&\hspace*{116pt}{}- \frac{1}{\gamma+iu + \zeta+ iz} \biggr) \hat {h}(z) \,dz \,du.\nonumber
\end{eqnarray}
To implement this transformation, one needs to find the Fourier inverse
of $x\mapsto e^{\zeta x}h(x)$ and evaluate the above integral
numerically, by truncating the domain of the integration
%
\begin{eqnarray}\label{eqgrRdef}
g_{r,R}(x) &=& \frac{1}{2\pi} \int
_{-r}^{r} e^{(\gamma+iu) x}
\int
_{-R}^{R} \biggl( \frac{\psi'(\gamma+iu)}{\psi(\gamma+iu) -
\psi(-\zeta-iz)}
\nonumber\\[-8pt]\\[-8pt]\nonumber
&&\hspace*{127pt}{}-
\frac{1}{\gamma+iu + \zeta+ iz} \biggr) \hat {h}(z) \,dz \,du,\hspace*{-10pt}
\end{eqnarray}
and, then, consider $r,R\rightarrow\infty$.
The next corollary provides the rate of convergence of $g_{r,R}$ to
$g$. Its proof is given in Appendix~\ref{appB}.

\begin{cor}\label{cormainrate}
Let Assumptions~\ref{ass1} and~\ref{ass2} hold, and let $g_{r,R}$
be given by (\ref{eqgrRdef}), with an arbitrary $\gamma> R_1$,
where $R_1$ is the constant appearing in Lemma~\ref{lepsiprimeUpsilonest}. Then, there exists a constant $c>0$,
independent of $\hat{h}$, $r$, $R$ and $x$, such that, for all large
enough $r,R>0$ and all $x>0$, we have
%
\begin{equation}
\label{eqmainrategrR} \bigl\llvert g(x) - g_{r,R}(x)\bigr\rrvert \leq c
\frac{e^{\gamma x}}{x} \biggl(\frac{1}{r}\llVert \hat{h}\rrVert
_{\mathbb{L}^1(\RR)} + \int_{\llvert  z\rrvert  >(r/2)\wedge R} \bigl\llvert \hat{h}(z)\bigr
\rrvert \,dz \biggr).
\end{equation}
\end{cor}


\section{Examples and implementation}\label{seexample}
Consider a L\'evy process $X$ given by the sum of a scaled Brownian
motion and a negative Gamma process. In other words,
%
\begin{equation}
\label{eqXBMGamma} X_t = \sigma B_t - \Gamma_t,
\end{equation}
where $B$ is a standard Brownian motion, and $\Gamma$ is a Gamma
process with parameters $\alpha>0$ and $\beta>0$; cf. \cite
{Bertoin,s99}. In this case, the characteristic triplet of
$X$, as defined in (\ref{eqLevyTriplet}), is given by
\[
\biggl(\mu=-\frac{\beta}{\alpha}, \sigma=\sigma, \Pi (dx)=\beta\frac{e^{-\alpha\llvert  x\rrvert  }}{\llvert  x\rrvert  }
\bone_{(-\infty,0)}(x)\,dx \biggr), %
\]
and in particular,
\[
\psi(\lambda) = \frac{1}{2}\sigma^2 \lambda^2 -
\beta\log \biggl(1 + \frac{\lambda}{\alpha} \biggr). %
\]

\begin{figure}

\includegraphics{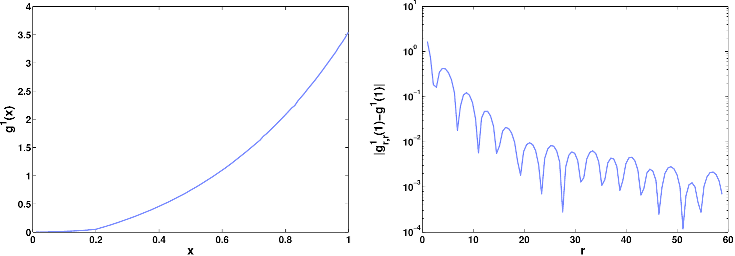}

\caption{On the left: the values of $g^1(x)$, defined in (\protect
\ref{eqgBMGamma}), for various $x>0$. On the right: the convergence
rate of the truncated integrals $g^1_{r,r}(1)$, defined in (\protect
\ref{eqgrRdef}), to the value of $g^1(1)$ (on a logarithmic scale).
The parameter values are: $K=-0.2$, $\alpha=\beta=\sigma=1$, $\zeta
= 0.9$, $\gamma=4$.}
\label{fig1}
\end{figure}

Consider an up-and-out put option, with maturity $T$, strike $K<0$ and
barrier $0$, written on the underlying process $S$, which has the
following payoff:
\[
(K - S_T )^+ \bone_{ \{\sup_{t\in[0,T]}S_t<0
\}}. %
\]
Here, for simplicity, we assume that the underlying process $S$ can
take negative values.
Assume that the risk neutral evolution of the underlying is given by
$S_t = x + X_t$, with some $x\in(-\infty,0)$ and with $X$ given by
(\ref{eqXBMGamma}).
Assume that we need to find a static hedging strategy for this barrier
option using the European-type options.
Then, following the algorithm presented in Sections~\ref{subseSRP}
and~\ref{subseWRP}, we need to construct the weak symmetry mapping
$\mathbf{W}^+$ and apply it to the hockey-stick function
$h^1(x)=(K-x)^+$. Note that Assumption~\ref{ass1} is always
satisfied for the process $X$, Assumption~\ref{ass2} holds with any
$\zeta\in(0,\alpha)$ and (\ref{eqtailcond}) is satisfied. Indeed,
it is easy to see that the inverse Fourier transform of the function
$x\mapsto e^{\zeta x} h(x)$ is given by
\[
\hat{h}^1(z) = \frac{\exp(K(\zeta+ iz))}{2\pi(\zeta+ iz)^2}. %
\]
Therefore, we apply the last assertion of Theorem~\ref{thmain} and
make use of (\ref{equopshformalpf}) to conclude that, at all times
up until and including the first time when $S$ hits $0$, the price of
up-and-out put option coincides with the price of a European-type
option, which has maturity $T$ and the following payoff function:
\[
(K- S_T )^+ - g_1 (S_T ), %
\]
with
%
\begin{eqnarray}\label{eqgBMGamma}
g^1(x) &=& \mathbf{W}^+ h^1(x) \nonumber
\\
&=& \frac{1}{4\pi^2 i}
\int_{\gamma-i\infty}^{\gamma+i\infty} e^{\lambda x} \int
_{\RR} \biggl(
\bigl(\sigma^2\lambda- \beta/(\lambda+\alpha)\bigr)\nonumber
\\
&&\hspace*{97pt}{}\big/
\biggl(\frac{\sigma^2\lambda^2}{2} - \beta\log \biggl(1+\frac{\lambda
}{\alpha} \biggr)
\nonumber\\[-8pt]\\[-8pt]\nonumber
&&\hspace*{108pt}{} - \frac{\sigma^2(\zeta+iz)^2}{2}+ \beta\log
 \biggl(1-\frac{\zeta+iz}{\alpha} \biggr)\biggr)
\\
&&\hspace*{218pt}{}
  -\frac{1}{\lambda+ \zeta
+ iz} \biggr)\nonumber
\\
&&\hspace*{60pt}{}\times \frac{e^{K(\zeta+ iz)}}{(\zeta+ iz)^2} \,dz \,d\lambda, \nonumber
\end{eqnarray}
for $x>0$.
Thus, in order to offset the risks associated with holding an
up-and-out put option (i.e., hedge the barrier option), one needs to
sell the\vadjust{\goodbreak} European-type option with the above payoff, and buy it back
(at a zero price) if and when the underlying hits $0$.
The results of numerical integration are presented in Figure~\ref{fig1}.
Notice that, in the present case, Corollary~\ref
{cormainrate} and the asymptotic relation $\llvert  \hat{h}^1(z)\rrvert  \sim\operatorname{const}\cdot\, z^{-2}$ imply that it is optimal to
approximate the double integral in (\ref{eqgBMGamma}) by\vspace*{1pt}
integrating over squares in the $\Im(\lambda)\times z$ domain. The
convergence rate of $g^1_{r,r}(1)$ [cf. (\ref{eqgrRdef})], as
$r\rightarrow\infty$, is shown on the right-hand side of Figure $\ref
{fig1}$ (on a logarithmic scale). Notice that the convergence rate
seems to be polynomial (rather than exponential), as predicted by
Corollary~\ref{cormainrate}. The numerical integration over a finite
domain, required to compute $g^1_{r,r}$, is performed via the MatLab
function \texttt{quad2d}. For the set of parameters used to generate
Figure~\ref{fig1}, the CPU time required to compute $g^1(1)_{r,r}$,
with $r=60$, is $1.27$ seconds.\footnote{All computations are
performed on a standard laptop, 1.8 GHz Intel Core i5, 4~GB RAM.}
%

\begin{figure}

\includegraphics{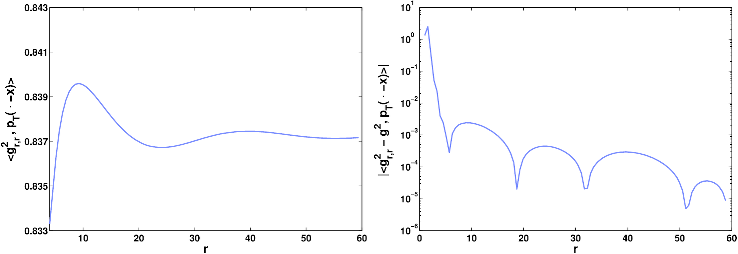}

\caption{On the left: the values of $\langle g^2_{r,r},p_T(\cdot
-x)\rangle$, with $g^2_{r,r}$ given by (\protect\ref
{eqgBMGamma}), for various $r>0$. On the right: the convergence rate
of $\langle g^2_{r,r},p_T(\cdot-x)\rangle$ to $\langle g^2,p_T(\cdot
-x)\rangle$, as $r\rightarrow\infty$ (on a logarithmic scale). The
parameter values are $x=0.1$, $T=1$, $K=-0.2$, $\alpha=\beta=\sigma
=1$, $\zeta= 0.9$, $\gamma=4$.}
\label{fig2}
\end{figure}

\begin{figure}

\includegraphics{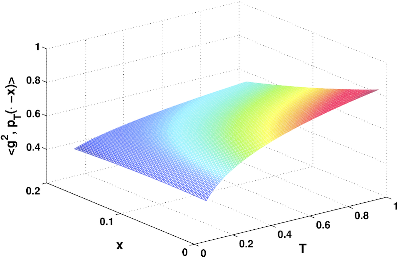}

\caption{The values of $\langle g^2_{r,r},p_T(\cdot-x)\rangle$, for
$x\in[0,0.2]$ and $T\in[0.1,1]$. The parameter values are: $r=60$,
$K=-0.2$, $\alpha=\beta=\sigma=1$, $\zeta= 0.9$, $\gamma=4$.}
\label{fig3}
\end{figure}

Next, we illustrate how the weak reflection principle can be used for
numerical computation of the joint marginal distribution of a
stochastic process and its running maximum. Assume that $X$ is given by
(\ref{eqXBMGamma}), and we need to approximate numerically the
value of
%
\begin{equation}
\label{eqBMGammajointlawprob} \PP \Bigl( X_T\leq K+x, \sup_{t\in[0,T]}
X_t \geq x \Bigr),
\end{equation}
with some $T>0$, $K<0$ and $x\geq0$.
The first assertion of Theorem~\ref{thmain} implies that
%
\begin{eqnarray}\label{eqBMGammajointlaw}
&& \PP \Bigl( X_T\leq K+x, \sup_{t\in[0,T]}
X_t \geq x \Bigr)\nonumber
\\
&&\qquad  = \EE h^2(X_T-x) - \EE
\bigl( h^2(X_T-x) \bone_{ \{ \sup_{t\in
[0,T]} (X_t-x) < 0 \}} \bigr)
\\
&&\qquad =\bigl\langle g^2, p_T(\cdot- x) \bigr\rangle,\nonumber
\end{eqnarray}
where $h^2=\bone_{(-\infty,K]}$ is the indicator function, $p_T$ is
the density of $X_T$ and $g^2=\mathbf{W}^+ h^2$ is the weak symmetry
image of $h^2$,
%
\begin{eqnarray}\label{eqBMGammajointlaw2}
&& \bigl\langle g^2, p_T(\cdot- x) \bigr\rangle\nonumber
\\
&&\qquad =
\frac{1}{2\pi i} \int_{\gamma-i\infty}^{\gamma+i\infty} \int
_{\RR} e^{\lambda y} p_T(y-x)\,dy \int
_{\RR} \biggl( \frac{\psi'(\lambda)}{\psi(\lambda) - \psi
(-\zeta-iz)} \nonumber
\\
&&\hspace*{228pt}{}- \frac{1}{\lambda+ \zeta+ iz} \biggr)
\hat{h}^2(z) \,dz \,d\lambda \nonumber
\\
&&\qquad  = \frac{1}{4\pi^2 i} \int_{\gamma-i\infty}^{\gamma+i\infty}
\exp \biggl(\lambda x + T \biggl(\frac{1}{2}\sigma^2
\lambda^2 - \beta\log (1 + \lambda/\alpha ) \biggr) \biggr)
\nonumber
\\
&&\hspace*{58pt}\quad\qquad{}\times
\int_{\RR} \biggl(
\bigl(\sigma^2\lambda- \beta/(\lambda+\alpha
)\bigr)
\\
&&\hspace*{123pt}{}\big/
\biggl(\frac{\sigma^2\lambda^2}{2} - \beta\log \biggl(1+\frac{\lambda
}{\alpha} \biggr) \nonumber
\\
&&\hspace*{134pt}{} - \frac{\sigma^2(\zeta+iz)^2}{2}+ \beta\log
 \biggl(1-\frac{\zeta+iz}{\alpha} \biggr)\biggr)\nonumber
 \\
 &&\hspace*{244pt}{}  -\frac{1}{\lambda+ \zeta
+ iz}
\biggr)\nonumber
\\
&&\hspace*{91pt}{}\times  \frac{e^{K(\zeta+ iz)}}{(\zeta+ iz)} \,dz \,d\lambda, \nonumber
\end{eqnarray}
and the associated $\hat{h}^2$ is given by
\[
\hat{h}^2(z) = -\frac{\exp(K(\zeta+ iz))}{2\pi(\zeta+ iz)}. %
\]
Notice that, strictly speaking, Theorem~\ref{thmain} yields (\ref
{eqBMGammajointlaw}) only for $x=0$.
To show that (\ref{eqBMGammajointlaw}) holds for all $x\geq0$, we
use $g^2_r$ in place of $g^2$ and repeat (\ref{equopshformalpf}),
to estimate the absolute value of the difference between the left and
the right-hand sides of (\ref{eqBMGammajointlaw}) by
\[
\sup_{t\in[0,T]} \bigl\llvert \bigl\langle g^2,
p_t\bigr\rangle- \EE g^2_r
(X_{t} )\bigr\rrvert = \sup_{t\in[0,T]} \bigl\llvert
\bigl\langle g^2-g^2_r, p_t\bigr
\rangle\bigr\rrvert. %
\]
Then we make use of (\ref{eqmainrategrweakpt}), to pass to the
limit as $r\rightarrow\infty$ and obtain (\ref{eqBMGammajointlaw}).
Figure~\ref{fig2} shows the convergence rate of the numerical
approximation of the right-hand side of (\ref
{eqBMGammajointlaw2}). As before, we used the MatLab function
\texttt{quad2d} to evaluate the integral in (\ref
{eqBMGammajointlaw2}) numerically, with the parameters' values
specified in Figure~\ref{fig2}. The CPU time required to compute
$\langle g^2_{r,r},p_T(\cdot-x)\rangle$, with $r=60$, $x=0.1$ and
$T=1$, is $1.69$ seconds.
%
It is worth mentioning that, as discussed in Section~\ref{subseWRP},
the computation of the joint marginal distribution of a process and its
running maximum is \emph{not} the main application of the weak
reflection principle. In fact, our method produces more than just the
value of the expectation of a function of $X_T$ and $\sup_{t\in
[0,T]}X_t$: it allows to express this value via the expectation of a
function of $X_T$ alone. The latter amounts to solving an inverse
problem, as opposed to the direct problem of computing the expectation.
Thus the use of the weak reflection principle for computing the joint
probabilities (as opposed to static hedging) does not fully utilize~the
power of the method. As a result, in general, our method may not
outperform the existing algorithms for the computation of the joint
law, which are based on the Wiener--Hopf factorization; cf. \cite
{s99,Kuznetsov3,Kuznetsov2,Kuznetsov1} and
references therein. Recall that the Wiener--Hopf factorization allows
one to compute the joint probability (\ref{eqBMGammajointlawprob})
by evaluating an integral over a vertical line on a complex plane, very
much like the integral with respect to $\lambda$ in (\ref
{eqBMGammajointlaw2}). The value of the integrand at each point is,
in turn, computed by a Fourier inversion applied to the Wiener--Hopf
factors; see, for example, equation (18) in \cite{Kuznetsov2}.\hskip.2pt\footnote{As described in \cite{Kuznetsov3}, for certain
families of L\'evy processes, the required Fourier inversion can be
reduced to a series expansion, which is more computationally
efficient.} Thus, in general, the computation of the integral in (\ref
{eqBMGammajointlaw2}) seems to have the same order of complexity as
the classical method [in addition, (\ref{eqBMGammajointlaw})--(\ref
{eqBMGammajointlaw2}) only apply to spectrally-negative L\'evy
processes, with nontrivial Brownian component]. However, in some
cases, the algorithm described by \mbox{(\ref{eqBMGammajointlaw})--(\ref{eqBMGammajointlaw2})}
may be more efficient. Namely, it is
advantageous to use the above method if the joint probability (\ref
{eqBMGammajointlawprob}) needs to be computed for multiple $x$ and
$T$ (which corresponds to varying the initial condition and the time
horizon). Notice that the inner integral on the right-hand side of
(\ref{eqBMGammajointlaw2}) is independent of $x$ and~$T$. Then the
computational complexity can be reduced by reusing, for different $x$
and $T$, the same values of the inner integral (as a function of
$\lambda$, computed on a given grid or via basis expansion), in the
integration with respect to $\lambda$.\footnote{Of course, this only
works if the $x$ and $T$ very over a reasonably small range of values,
so that there is no need to change the precision with which the inner
integral in (\ref{eqBMGammajointlaw2}) is computed. Also, by
changing the order of integration, one can use the same idea to compute
the right-hand side of (\ref{eqBMGammajointlaw2}) for multiple $K$
(with $x$ and $T$ fixed) in a more efficient way.} Figure~\ref{fig3}
shows the numerical approximation of $\langle g^2, p_T(\cdot- x)
\rangle$, given by the right-hand\vadjust{\goodbreak} side of~(\ref{eqBMGammajointlaw2}), for 10,000 different pairs $(x,T)$. In
this computation, we approximate the integral with respect to $\lambda
$, in (\ref{eqBMGammajointlaw2}), by a simple Riemann sum (with
the uniform partition of diameter $0.6$) and reuse the same values of
the inner integral in~(\ref{eqBMGammajointlaw2}), for different
$x$ and $T$ (the inner integral is computed using the MatLab function
\texttt{quad}). As a result, the total CPU time required to finish all
the computations is only $20.46$ seconds [compared to $1.69$ seconds
required to evaluate (\ref{eqBMGammajointlaw2}) for a single pair $(x,T)$].
Another advantage of (\ref{eqBMGammajointlaw})--(\ref
{eqBMGammajointlaw2}) is that these formulas enable a
straightforward computation of the derivatives of the joint probability
(\ref{eqBMGammajointlawprob}) with respect to $x$ and $T$ (which
provide sensitivities with respect to the initial condition and the
time horizon). Indeed, it follows from (\ref{eqfourdensasymp}) and
(\ref{eqcor2est1}) that the absolute value of the integrand on the
right-hand side of (\ref{eqBMGammajointlaw2}) decays faster than
any exponential, as a function of $\Im(\lambda)$. Therefore, the same
conclusion holds for any derivative of the integrand with respect to
$x$ and $T$. Thus we can differentiate with respect to $x$ and $T$
inside the integral in (\ref{eqBMGammajointlaw2}) an arbitrary
number of times.

%

Finally, we illustrate the PIDE interpretation of the weak symmetry.
Recall that every L\'evy process has an infinitesimal generator
$\mathcal{L}$ associated with it; cf. \cite{s99}. This generator can
be viewed as a pseudo-differential operator, which acts on all
infinitely smooth functions $\phi$ with compact support as follows:
\[
\mathcal{L}\phi(x) = \frac{\sigma^2}{2} \,\partial^2_x
\phi(x) + \mu\,\partial_x \phi(x) + \int_{\RR}
\bigl( \phi(x+z) - \phi(x) - z\,\partial_x\phi(x) \bigr) \Pi(dz).
\]
In the present case, we obtain
\begin{eqnarray*}
\mathcal{L}\phi(x) &=& \frac{\sigma^2}{2} \,\partial^2_x
\phi(x) -\frac{\beta}{\alpha} \,\partial_x \phi(x)
\\
&&{}- \beta\int
_{-\infty}^0 \bigl( \phi(x+z) - \phi(x) - z\,\partial
_x\phi(x) \bigr) \frac{e^{\alpha z}}{z}\,dz. %
\end{eqnarray*}
Let us introduce
\[
u(x,t) = \EE h (x+X_t ). %
\]
Then the Feynman--Kac formula for the L\'evy process $X$ (see, e.g.,
\cite{Lions,Cont,Hilber} and \cite{Eberlein})
implies that $u$ is a solution of the following initial value problem:
%
\begin{equation}
\label{eqinitValdef} \cases{ \partial_t u - \mathcal{L} u = 0, &\quad $x
\in\RR$, $t>0$,
\cr
u(x,0) = h(x).}
\end{equation}
The exact definition of a solution to the above problem depends on the
regularity assumptions on function $h$ and is discussed in the
aforementioned references.
The existence of the weak symmetry mapping $\mathbf{W}^+$ for the
process $X$ implies that, for any admissible $h$, with support in
$(-\infty,0)$, we can find a function $g = \mathbf{W}^+ h$, with
support in $[0,\infty)$, such that the solution to (\ref
{eqinitValdef}), with $h$ replaced by $g$, coincides with the
original solution at $x=0$, for all $t>0$. In particular, this means
that we can modify any initial condition for $x\in[0,\infty)$ to
ensure that $u(0,t)=0$, for all $t>0$. This, in turn, implies that we
can reduce an initial-boundary value problem to an initial-value problem.
Namely, any solution to the PIDE in (\ref{eqinitValdef}), defined
for $(x,t)\in(-\infty,0)\times(0,\infty)$, with zero boundary
condition at $x=0$ and with initial condition $h$, can be represented
as a solution to the initial value problem~(\ref{eqinitValdef}),
with $h$ replaced by $h-g$.
Notice that equations of the form (\ref{eqinitValdef}) may be
interesting on their own. For example, it is discussed in \cite
{Carmona} and \cite{Baeumer} that the relativistic Schr\"odinger
equation can be reduced to a PIDE associated with the L\'evy process
known as the normal inverse Gaussian (NIG) process. Unfortunately,
Assumption~\ref{ass1} excludes NIG from the scope of the present
work, since NIG is a pure jump process with both positive and negative
jumps.\footnote{Note that we focused on the L\'evy processes with only
negative jumps to ensure that the process does not jump across the
barrier. Indeed, this is necessary for the weak reflection principle to
hold, as we need to stop the process precisely at the barrier. However,
in the PIDE interpretation, we do not need to use all the steps of the
weak reflection principle. Namely, we do not need to stop the process
at the barrier, but rather, we only need to construct the weak symmetry
mapping, assuming the process starts from the barrier. Hence, in the
context of PIDEs, it makes sense to consider processes with both
positive and negative jumps.} However, we believe that future research
will extend the results presented herein, to construct weak symmetry
mappings for NIG and other important L\'evy processes with two-sided
jumps and/or missing Brownian component.

\section{Summary and extensions}\label{sesummary}

We have presented a new mathematical technique, which we christened the
\emph{weak reflection principle} and which is an extension of the
well-known reflection principle for Brownian motion. This new form of
reflection principle is obtained by weakening the notion of symmetry
that is required for the classical reflection principle to hold. More
precisely, our method is based on the notion of \emph{weak symmetry}.
We started by reviewing the existing results which provide an explicit
integral representation of the weak symmetry mapping for any
time-homogeneous diffusion process on a real line (subject to some
regularity conditions). Finally, for the most of the paper we focused
on constructing the weak symmetry mapping for spectrally negative L\'
evy processes, thus extending the weak reflection principle to this new
class of stochastic processes.

The weak reflection principle provides solutions to various problems
for which the classical reflection principle can be used, even when the
underlying process is not a Brownian motion and does not possess any
strong symmetries. 
In particular, the weak reflection principle is a perfect tool for
constructing the exact static hedging strategies of barrier options (in
fact, this problem motivated the development of the method in the first
place). Another application of this method is the computation of the
joint distribution of a process and its running maximum (minimum). Of
course, while this problem is quite relevant for diffusions, in the
case of L\'evy processes, there exist several alternative computational
methods, based on the Wiener--Hopf factorization; cf. \cite{s99,Kuznetsov3,Kuznetsov2,Kuznetsov1}.
Nevertheless, as shown in Section~\ref{seexample}, there are cases
when it is advantageous to use this method.
Finally, the weak symmetry mapping allows us to solve an inverse
problem for the parabolic partial PIDE associated with a L\'evy
process. Namely, using the weak symmetry, we can modify the initial
condition of the PIDE on one half line only, so that its solution
remains constant at $x=0$, for all times. This, in particular, allows
us to represent the solution of a PIDE with initial and boundary
conditions via the solution of the same PIDE with initial condition only.

It is also worth mentioning that the technical Lemmas~\ref
{lepsiest1}--\ref{lepsiprimeUpsilonest} and the resulting
Corollary~\ref{coranExt} describe a domain on which the resolvent
function of a L\'evy process is well defined. This domain is rather
large, and in particular, the real parts of its elements are unbounded
from below. Thus our results provide a nontrivial estimate of the
spectrum of the integro-differential operator associated with any
admissible spectrally-negative L\'evy process. Recall that this
operator is nonlocal and nonsymmetric, which makes it very hard to
describe its spectrum using the general theory; see also Remark~\ref{remrem1} for a description of the associated difficulties.

To date, the weak symmetry has only been established for diffusion
processes and L\'evy processes with one-sided jumps. However, we
conjecture that these results can be extended to a larger class of
time-homogeneous Markov processes---possibly all jump-diffusions
satisfying some regularity conditions. Such an extension would allow us
to solve the aforementioned problems for a larger class of stochastic processes.
In particular, as discussed in Section~\ref{seexample}, including the
NIG process in the scope of our analysis would establish a connection
with the relativistic Schr\"odinger equation which is important in Physics.
Another possible extension is related to the domain with respect to
which the weak symmetry is defined. Notice that, in the present case,
we split the real line into two half lines and study the weak symmetry
of the process with respect to the (unique) boundary point. It is
interesting to extend these results to the case of a compact interval,
whose boundary consists of two points (assuming the underlying Markov
process does not jump across the boundary points). In financial
mathematics, this problem would correspond to the static hedging of
double barrier options. More generally, one can investigate domains in
higher dimension and try to establish the weak symmetry with respect to
their boundaries.

\begin{appendix}
\section{}\vspace*{-12pt}\label{appA}

\begin{pf*}{Proof of Lemma~\ref{leIm1}}
\begin{eqnarray*}
&& \Im \bigl(\psi(v + iu) - \psi(-\zeta-iz) \bigr) %
\\
&&\qquad
= \mu(u+z) + (v u - z\zeta ) \sigma^2 + \int_{-\infty}^0
\bigl(e^{v x} \sin(ux) - u x \bigr) \Pi(dx)
\\
&&\quad\qquad{} + \int_{-\infty}^0
\bigl(e^{-\zeta x} \sin(zx) - zx \bigr) \Pi(dx). %
\end{eqnarray*}
Notice that
\begin{eqnarray*}
&& \int_{-\infty}^0 \bigl(e^{v x} \sin(ux) - u
x \bigr) \Pi(dx)
\\
&&\qquad = v u \int_{-\infty}^0
\frac{e^{v x} - 1}{v x} \frac{\sin(ux)}{ux} x^2 \Pi(dx) + u^2
\int_{-\infty}^0 \frac{\sin(ux) - ux}{u^2x^2} x^2
\Pi(dx),
\\
&& \int_{-\infty}^0 \bigl(e^{-\zeta x} \sin(zx) -
zx \bigr) \Pi(dx)
\\
&&\qquad = z\zeta\int_{-\infty}^0
\frac{1 - e^{\zeta x}}{\zeta x} \frac
{\sin(zx)}{z x} x^2 e^{-\zeta x} \Pi(dx) +
z^2 \int_{-\infty}^0 \frac{ \sin(zx) - zx }{z^2 x^2}
x^2\Pi(dx).
\end{eqnarray*}
Since $x^2$ and $x^2 e^{-\zeta x}$ are integrable with respect to $\Pi
(dx)$, due to the dominated convergence theorem, the above integrals
are absolutely bounded and vanish as the corresponding functions of
$(u,v,z)$ go to infinity.
\end{pf*}

\begin{pf*}{Proof of Lemma~\ref{leIm2}}
We follow the proof of Lemma~\ref{leIm1}, except that at the end,
we apply the following additional estimate:
\begin{eqnarray*}
&& u^2\int_{-\infty}^0 \frac{\sin(ux) - ux}{u^2x^2}
x^2 \Pi(dx) + z^2 \int_{-\infty}^0
\frac{ \sin(zx) - zx }{z^2 x^2} x^2\Pi(dx) %
\\
&&\qquad = \int_{-\infty}^0 \bigl( 2\sin\bigl((z+u)x/2\bigr)
\cos\bigl((z-u)x/2\bigr) - (z+u)x \bigr) \Pi(dx) %
\\
&&\qquad = \bigl(z^2-u^2\bigr)\int_{-\infty}^0
2\frac{\sin((z+u)x/2)}{(z+u) x} \frac
{\cos((z-u)x/2) - 1}{(z-u)x} x^2 \Pi(dx) %
\\
&&\quad\qquad{}+ (z+u)^2 \int_{-\infty}^0
\frac{2\sin((z+u)x/2) -
(z+u)x}{(z+u)^2x^2} x^2 \Pi(dx). %
\end{eqnarray*}
Applying the dominated convergence theorem, we complete the proof of
the lemma.
\end{pf*}

\begin{pf*}{Proof of Lemma~\ref{leRe}}
\begin{eqnarray*}
&& \Re \bigl(\psi(v + iu) - \psi(-\zeta-iz) \bigr)
\\
&&\qquad = \mu(v+\zeta) +
\bigl(v^2 - \zeta^2 + z^2 - u^2
\bigr)\sigma^2/2
\\
&&\quad\qquad{}
+ \int_{-\infty}^0 \bigl(e^{v x} \cos(ux) -
1 -v x \bigr) \Pi(dx)
\\
&&\quad\qquad{} - \int_{-\infty}^0
\bigl(e^{-\zeta x} \cos(zx) - 1 + \zeta x \bigr) \Pi(dx). %
\end{eqnarray*}
Notice that
\begin{eqnarray*}
&& \int_{-\infty}^0 \bigl(e^{v x} \cos(ux) - 1
- v x \bigr) \Pi(dx) - \int_{-\infty}^0
\bigl(e^{-\zeta x} \cos(zx) - 1 + \zeta x \bigr) \Pi(dx) %
\\
&&\qquad = \int_{-\infty}^0 e^{v x} \bigl(\cos(ux) -
1\bigr) \Pi(dx) + \int_{-\infty}^0
\bigl(e^{v x} - 1 - v x \bigr) \Pi(dx) %
\\
&&\quad\qquad{}- \int_{-\infty}^0 e^{-\zeta x} \bigl(\cos(zx) -
1\bigr) \Pi(dx) - \int_{-\infty}^0
\bigl(e^{-\zeta x} - 1 + \zeta x \bigr) \Pi(dx) %
\\
&&\qquad = \int_{-\infty}^0 \bigl(\cos(ux) - 1\bigr) \Pi(dx)
+ uv \int_{-\infty}^0 \frac{e^{v x}-1}{v x}
\frac{\cos(ux) -
1}{ux} x^2 \Pi(dx)
\\
&&\quad\qquad{}+v^2 \int
_{-\infty}^0 \frac{e^{v x} - 1 - v x }{v^2 x^2} x^2
\Pi(dx)
- \int_{-\infty}^0 \bigl(\cos(zx) - 1\bigr) \Pi(dx)
\\
&&\quad\qquad{}
- z\zeta\int_{-\infty}^0 \frac{1 - e^{\zeta x}}{\zeta x}
\frac
{\cos(zx) - 1}{zx} x^2 e^{-\zeta x} \Pi(dx)
\\
&&\quad\qquad{} - \zeta^2
\int_{-\infty}^0 \frac{1 - e^{\zeta x} + \zeta x e^{\zeta
x}}{\zeta^2 x^2} x^2
e^{-\zeta x} \Pi(dx) %
\\
&&\qquad = \int_{-\infty}^0 \bigl(\cos(ux) - \cos(zx)\bigr)
\Pi(dx) + \bigl\llvert \overline{\overline{o}}(uv)\bigr\rrvert + \bigl\llvert
\overline {\overline{o}}\bigl(v^2\bigr)\bigr\rrvert + \overline{
\overline{o}}(z) + \overline{\overline{o}}(1) %
\\
&&\qquad = 2\bigl(z^2-u^2\bigr) \int_{-\infty}^0
\frac{\sin((u+z)x/2)}{(z+u)x} \frac
{\sin((z-u)x/2)}{(z-u)x} x^2 \Pi(dx)
\\
&&\quad\qquad{} + \bigl\llvert
\overline{\overline{o}}(uv)\bigr\rrvert + \bigl\llvert \overline {\overline{o}}
\bigl(v^2\bigr)\bigr\rrvert + \overline{\overline{o}}(z) + \overline{
\overline{o}}(1) %
\\
&&\qquad = \overline{\overline{o}}\bigl(z^2-u^2\bigr) + \bigl
\llvert \overline{\overline{o}}(uv)\bigr\rrvert + \bigl\llvert \overline {
\overline{o}}\bigl(v^2\bigr)\bigr\rrvert + \overline{\overline{o}}(z)
+ \overline{\overline{o}}(1), %
\end{eqnarray*}
where we applied the dominated convergence theorem.
\end{pf*}

\begin{pf*}{Proof of Lemma~\ref{lepsiest1}}
Using Lemmas~\ref{leIm2} and~\ref{leRe}, we obtain
%
\begin{eqnarray}\label{eqpsiest11}
&& \bigl\llvert \psi(v + iu) - \psi(-\zeta-iz)\bigr\rrvert\nonumber
\\
&&\qquad  \geq \bigl\llvert \Im
\bigl(\psi(v + iu) - \psi(-\zeta-iz) \bigr)\bigr\rrvert \bone _{ \{\llvert  \llvert  u\rrvert   - \sqrt{z^2 + v^2} \rrvert   \leq\varepsilon\nonumber
v \}}
\\
&&\quad\qquad{} + \bigl\llvert \Re \bigl(\psi(v + iu) - \psi(-\zeta-iz) \bigr)\bigr
\rrvert \bone_{ \{\llvert  \llvert  u\rrvert   - \sqrt{z^2 + v^2} \rrvert   > \varepsilon
v \}}\nonumber
\\
&&\qquad = \bigl\llvert \mu(u+z) + (v u - z \zeta) \sigma^2
+ \overline{
\overline {o}}(v u)+ \overline{\overline{o}}(z)
\\
&&\hspace*{36pt}{} + \overline{\overline{o}}
\bigl(z^2-u^2\bigr)+ \operatorname{sign}(z+u) \bigl
\llvert \overline{\overline{o}} \bigl((z+u)^2 \bigr)\bigr\rrvert \bigr
\rrvert \bone_{ \{\llvert  \llvert  u\rrvert   - \sqrt{z^2 + v^2} \rrvert   \leq
\varepsilon v \}}\nonumber
\\
&&\quad\qquad{}
+ \bigl\llvert \mu(v+\zeta) + \bigl(v^2 + z^2 -
u^2 \bigr)\sigma^2/2 \nonumber
\\
&&\hspace*{48pt}{} + \overline{\overline{o}}
\bigl(z^2-u^2\bigr) + \bigl\llvert \overline{
\overline{o}}(uv)\bigr\rrvert + \bigl\llvert \overline {\overline{o}}
\bigl(v^2\bigr)\bigr\rrvert + \overline{\overline{o}}(z) + \overline{
\overline{o}}(1) \bigr\rrvert \bone_{ \{\llvert  \llvert  u\rrvert   - \sqrt{z^2 + v^2} \rrvert   > \varepsilon
v \}}.\hspace*{-10pt}\nonumber
\end{eqnarray}
Let us estimate the first term in the above:
\begin{eqnarray*}
&& \bigl\llvert \mu(u+z) + v u \sigma^2 - z \zeta\sigma^2 +
\overline {\overline{o}}(v u)+ \overline{\overline{o}}(z)
\\
&&\quad{}  + \overline{
\overline{o}}\bigl(z^2-u^2\bigr) + \operatorname{sign}(z+u)
\bigl\llvert \overline{\overline{o}} \bigl((z+u)^2 \bigr)\bigr\rrvert
\bigr\rrvert \bone_{ \{\llvert  \llvert  u\rrvert   - \sqrt{z^2 + v^2} \rrvert   \leq
\varepsilon v \}} %
\\
&&\qquad \geq \bigl( v \llvert u\rrvert \sigma^2 + \mu u + z\bigl(\mu- \zeta
\sigma^2\bigr) + \overline {\overline{o}}\bigl(v \bigl(\llvert z
\rrvert +v\bigr)\bigr)
\\
&&\hspace*{35pt}{} + \overline{\overline{o}}(z) + \overline{\overline{o}}
\bigl(v\llvert z\rrvert + v^2\bigr) + \operatorname{sign}(z+u) \bigl
\llvert \overline{\overline{o}} \bigl((z+u)^2 \bigr)\bigr\rrvert
\bigr)
\bone_{ \{\llvert  \llvert  u\rrvert   - \sqrt{z^2 + v^2} \rrvert   \leq
\varepsilon v, u\geq0 \}} %
\\
&&\quad\qquad{} + \bigl( v \llvert u\rrvert \sigma^2 - \mu u - z\bigl(\mu- \zeta
\sigma^2\bigr)
\\
&&\hspace*{48pt}{} + \overline {\overline{o}}\bigl(v \bigl(\llvert z
\rrvert +v\bigr)\bigr)
+ \overline{\overline{o}}(z) + \overline{\overline{o}}
\bigl(v\llvert z\rrvert + v^2\bigr)
\\
&&\hspace*{97pt}{} - \operatorname{sign}(z+u) \bigl\llvert \overline{\overline{o}} \bigl((z+u)^2 \bigr)\bigr\rrvert
\bigr)
\bone_{ \{\llvert  \llvert  u\rrvert   - \sqrt{z^2 + v^2} \rrvert   \leq
\varepsilon v, u<0 \}}
\\
&&\qquad \geq \bigl( c_3 v \bigl(\llvert z\rrvert +v\bigr) + \mu u + z\bigl(
\mu- \zeta\sigma^2\bigr) + \overline {\overline{o}}\bigl(v \bigl(
\llvert z\rrvert +v\bigr)\bigr)
\\
&&\hspace*{115pt}{}+ \overline{\overline{o}}(z) + \overline{
\overline{o}}\bigl(v\llvert z\rrvert + v^2\bigr) + \overline{\overline
{o}}\bigl(v^2\bigr) \bigr) \bone_{ \{\llvert  \llvert  u\rrvert   - \sqrt{z^2 + v^2} \rrvert   \leq
\varepsilon v, u\geq0 \}} %
\\
&&\qquad\quad{} + \bigl( c_4 v \bigl(\llvert z\rrvert +v\bigr) - \mu u - z\bigl(
\mu- \zeta\sigma^2\bigr)
\\
&&\hspace*{48pt}{} + \overline {\overline{o}}\bigl(v \bigl(
\llvert z\rrvert +v\bigr)\bigr)+ \overline{\overline{o}}(z) + \overline{
\overline{o}}\bigl(v\llvert z\rrvert + v^2\bigr) + \overline{\overline
{o}} \bigl(v^2 \bigr) \bigr) \bone_{ \{\llvert  \llvert  u\rrvert   - \sqrt{z^2 + v^2} \rrvert   \leq
\varepsilon v, u<0 \}} %
\\
&&\qquad \geq c_1 v \bigl(\llvert z\rrvert +v\bigr) \bone_{ \{\llvert  \llvert  u\rrvert   - \sqrt{z^2 + v^2} \rrvert
\leq\varepsilon v \}},
\end{eqnarray*}
where $c_i$'s are\vspace*{1pt} positive constants, and we assume that $\llvert  z\rrvert  $ and
$v>0$ are large enough. In the above, we use the fact that if $\llvert  \llvert  u\rrvert   - \sqrt{z^2 + v^2} \rrvert   \leq\varepsilon v$, then
\[
\frac{1-\varepsilon}{2} \bigl(v + \llvert z\rrvert \bigr) \leq\llvert u\rrvert \leq v
+ (1+\varepsilon) \llvert z\rrvert, %
\]
and if, in addition, $uz<0$, then
\[
\bigl\llvert \llvert u\rrvert - \llvert z\rrvert \bigr\rrvert \leq\varepsilon v
+ \frac{v^2}{\sqrt{z^2 +
v^2} + \llvert  z\rrvert  } \leq(1+\varepsilon) v. %
\]

Finally, we estimate the second term on the right-hand side of (\ref
{eqpsiest11}):
\begin{eqnarray*}
&& \bigl\llvert \mu(v+\zeta) + \bigl(v^2 + z^2 -
u^2 \bigr)\sigma^2/2 + \overline{\overline{o}}
\bigl(z^2-u^2\bigr) + \bigl\llvert \overline{
\overline{o}}(uv)\bigr\rrvert + \bigl\llvert \overline {\overline{o}}
\bigl(v^2\bigr)\bigr\rrvert + \overline{\overline{o}}(z) \bigr\rrvert
\\
&&\quad{}\times
\bone_{ \{\llvert  \llvert  u\rrvert   - \sqrt{z^2 + v^2} \rrvert   > \varepsilon
v \}} %
\\
&&\qquad \geq \bigl\llvert u^2 - \bigl(v^2 + z^2
\bigr) \bigr\rrvert
\\
&&\quad\qquad{}\times \biggl(\sigma^2/2 + \frac{\mu(v+z) + \overline{\overline{o}}(u^2-z^2)
+ \llvert  \overline{\overline{o}}(uv)\rrvert
+ \llvert  \overline{\overline{o}}(v^2)\rrvert
+ \overline{\overline{o}}(z)}{\llvert  u^2 - (v^2 + z^2) \rrvert  } \biggr)
\\
&&\quad\qquad{}\times
\bone_{ \{\llvert  \llvert  u\rrvert   - \sqrt{z^2 + v^2} \rrvert   > \varepsilon
v \}} %
\\
&&\qquad \geq c_2 \bigl\llvert u^2 - \bigl(v^2 +
z^2\bigr) \bigr\rrvert \bone_{ \{\llvert  \llvert  u\rrvert   - \sqrt{z^2 + v^2} \rrvert   > \varepsilon
v \}}, %
\end{eqnarray*}
which holds for all large enough $\llvert  z\rrvert  $ and $v>0$, with some positive
constants $c_i$.
In the above we make use of the fact that if $\llvert  \llvert  u\rrvert   - \sqrt{z^2 +
v^2} \rrvert   > \varepsilon v$, then
\begin{eqnarray*}
\llvert u\rrvert &\geq& (1-\varepsilon) v,
\\
\bigl\llvert u^2 -
\bigl(v^2 + z^2\bigr) \bigr\rrvert &\geq&\varepsilon v
\bigl(\llvert u\rrvert + \sqrt{z^2 + v^2} \bigr),
\\
\frac{\llvert  u^2-z^2\rrvert  }{\llvert  u^2 - (v^2 + z^2) \rrvert  } &\leq&\frac{v}{v^2-v} \bone_{ \{ \llvert  u^2-z^2\rrvert  \leq v  \}}
+ \biggl(1 +
\frac{v^2}{\varepsilon v  (\llvert  u\rrvert   + \sqrt{v^2 + z^2}
 )} \biggr) \bone_{ \{ \llvert  u^2-z^2\rrvert  > v  \}} %
\\
&\leq&\frac{1}{v-1} + \biggl(1 + \frac{1}{\varepsilon} \biggr)
\bone_{ \{
\llvert  u^2-z^2\rrvert  > v  \}}, %
\end{eqnarray*}
and hence
\[
\frac{\mu(v+z) + \overline{\overline{o}}(u^2-z^2)
+ \llvert  \overline{\overline{o}}(uv)\rrvert
+ \llvert  \overline{\overline{o}}(v^2)\rrvert
+ \overline{\overline{o}}(z)}{\llvert  u^2 - (v^2 + z^2) \rrvert  } %
\]
can be made arbitrarily small by choosing large enough $v$.
\end{pf*}

\begin{pf*}{Proof of Lemma~\ref{lepsiest2}}
First, we notice that since $\psi$ is analytic, we have $\llvert  \psi
(v-iu)\rrvert  =\llvert  \psi(v+iu)\rrvert  $, and hence it suffices to consider only $u>0$.

Using Lemma~\ref{leRe}, we obtain the following inequalities:
\begin{eqnarray*}
&& \bigl( \Re \bigl( \psi(v + iu) - \psi(-\zeta-iz) \bigr) \bigr)^2
\\
&&\qquad = \bigl(\mu(v+\zeta) + \bigl(v^2 + z^2 - u^2
\bigr)\sigma^2/2 + \overline{\overline{o}} \bigl(z^2 -
u^2 \bigr)
\\
&&\hspace*{98pt}{}+ \overline{\overline{o}}(uv) + \overline{\overline{o}}
\bigl(v^2 \bigr) + \overline{\overline{o}}(z) + \overline{\overline
{o}}(1) \bigr)^2 %
\\
&&\qquad\geq \bigl(v^2 + (v+\zeta)2\mu/\sigma^2 +
z^2 - u^2 \bigr)^2\sigma^4/4
\\
&&\quad\qquad{}+ \bigl(v^2 + z^2 - u^2 + (v+\zeta)2\mu/
\sigma^2 \bigr)
\\
&&\qquad\qquad{}\times  \bigl(\overline{\overline{o}} \bigl(z^2 -
u^2 \bigr) + \overline{\overline{o}}(uv) + \overline{\overline{o}}
\bigl(v^2 \bigr) + \overline{\overline{o}}(z) \bigr)
\bone_{ \{\llvert  u-\llvert  z\rrvert  \rrvert  \leq\varepsilon \}} %
\\
&&\quad\qquad{} + \bigl(v^2 + z^2 - u^2 + (v+\zeta)2\mu/
\sigma^2 \bigr)
\\
&&\qquad\qquad{}\times \bigl( \overline{\overline{o}} \bigl(z^2 -
u^2 \bigr) + \overline{\overline{o}}(uv) + \overline{\overline{o}}
\bigl(v^2 \bigr) + \overline{\overline{o}}(z) \bigr)
\bone_{ \{\llvert  u-\llvert  z\rrvert  \rrvert  >\varepsilon \}} 
\\
&&\qquad \geq \bigl(v^2 + (v+\zeta)2\mu/\sigma^2
\bigr)^2\sigma^4/4
\\
&&\quad\qquad{}+ \bigl(v^2 + (v+\zeta)2
\mu/\sigma^2 \bigr)
 \bigl(z^2 - u^2\bigr)
\sigma^4/2 + \bigl(z^2 - u^2
\bigr)^2\sigma^4/4 
\\
&&\quad\qquad{}+ \bigl(v^2 + u \bigr) \bigl(\overline{\overline{o}}
\bigl(z^2 - u^2 \bigr) + \overline{\overline{o}}(uv) +
\overline{\overline{o}} \bigl(v^2 \bigr) + \overline{\overline{o}}(z)
\bigr) \bone_{ \{\llvert  u-\llvert  z\rrvert  \rrvert  \leq\varepsilon \}} %
\\
&&\quad\qquad{} + \bigl(v^2 + z^2 - u^2 + (v+\zeta)2\mu/
\sigma^2 \bigr)
\\
&&\qquad\qquad{}\times  \bigl( \overline{\overline{o}} \bigl(z^2 -
u^2 \bigr) + \overline{\overline{o}}(uv) + \overline{\overline{o}}
\bigl(v^2 \bigr) + \overline{\overline{o}}(z) \bigr)
\bone_{ \{\llvert  u-\llvert  z\rrvert  \rrvert  >\varepsilon \}} %
\\
&&\qquad = \bigl(v^2 + (v+\zeta)2\mu/\sigma^2
\bigr)^2\sigma^4/4 + \bigl(z^2 -
u^2 \bigr)^2\sigma^4/4
\\
&&\quad\qquad{}+
\bigl(v^2 + (v+\zeta)2\mu/\sigma^2 \bigr) z^2
\sigma^4/2 %
- \bigl(v^2 + (v+\zeta)2\mu/\sigma^2 \bigr)
u^2\sigma^4/2 %
\\
&&\quad\qquad{}+ \bigl(v^2 + u \bigr) \bigl(u\overline{\overline{o}} (v ) +
\overline{\overline{o}} \bigl(v^2 \bigr) + \overline{\overline {o}}(u)
\bigr) \bone_{ \{\llvert  u-\llvert  z\rrvert  \rrvert  \leq\varepsilon \}} %
\\
&&\quad\qquad{}+ \bigl(v^2 + \bigl\llvert z^2 - u^2\bigr
\rrvert \bigr) \bigl( \overline{\overline {o}} \bigl(z^2 -
u^2 \bigr) + \overline{\overline{o}}(uv) + \overline{\overline{o}}
\bigl(v^2 \bigr) + \overline{\overline{o}}(z) \bigr)
\bone_{ \{\llvert  u-\llvert  z\rrvert  \rrvert  >\varepsilon \}} %
\\
&&\qquad \geq v^4 \bigl(1-\varepsilon'\bigr)
\sigma^4/4 + \bigl(z^2 - u^2
\bigr)^2\sigma^4/4
+ v^2 z^2
\bigl(1-\varepsilon'\bigr)\sigma^4/2
\\
&&\quad\qquad{} - u^2
v^2 \bigl(1+\varepsilon'\bigr)\sigma^4/2
+ u^2\overline{\overline{o}} (v ) + \overline{\overline{o}}
\bigl(u^2 \bigr) + u\overline{\overline{o}} \bigl(v^3
\bigr) + v^2\overline{\overline{o}} (u ) + \overline{\overline{o}}
\bigl(v^4 \bigr) %
\\
&&\quad\qquad{}+ \bigl(v^2 + \bigl\llvert z^2 - u^2\bigr
\rrvert \bigr) \bigl( \overline{\overline {o}} \bigl(z^2 -
u^2 \bigr) + \overline{\overline{o}}(uv) + \overline{\overline{o}}
\bigl(v^2 \bigr) + \overline{\overline{o}}(z) \bigr)
\bone_{ \{\llvert  u-\llvert  z\rrvert  \rrvert  >\varepsilon \}}, %
\end{eqnarray*}
which hold for any (fixed) $\varepsilon,\varepsilon'>0$ and all large
enough $u$, $v$ and $\llvert  z\rrvert  $.

Next, using Lemmas~\ref{leIm1},~\ref{leIm2}, we obtain the
following inequalities:
\begin{eqnarray*}
&& \bigl( \Im \bigl( \psi(v + iu) - \psi(-\zeta-iz) \bigr) \bigr)^2
\\
&&\qquad = \bigl(\mu(u+z) + uv\sigma^2 - z\zeta\sigma^2 +\overline {\overline{o}}(uv) + \overline{\overline{o}}(z)
\\
&&\hspace*{97pt}{}  + \overline{
\overline{o}}\bigl(z^2-u^2\bigr) + \overline{\overline{o}}
\bigl((z+u)^2 \bigr) \bigr)^2 \bone_{ \{z\leq0 \}}
\\
&&\quad\qquad{}+ \bigl(\mu(u+z) + uv\sigma^2 - z\zeta\sigma^2 +
\overline {\overline{o}}(v u) + \overline{\overline{o}}(z) + \bigl\llvert
\overline{\overline{o}}\bigl(u^2\bigr) \bigr\rrvert + \bigl\llvert
\overline {\overline{o}}\bigl(z^2\bigr) \bigr\rrvert
\bigr)^2 \bone_{ \{z>0 \}} %
\\
&&\qquad \geq \bigl( uv\sigma^2 + \mu(u+z) - z\zeta\sigma^2 +
\overline {\overline{o}}(v u) + \overline{\overline{o}}(z) \bigr)^2
\\
&&\quad\qquad{}+ \bigl(\mu(u+z) + uv \sigma^2 - z\zeta\sigma^2 +
\overline {\overline{o}}(v u) + \overline{\overline{o}}(z) \bigr)
\\
&&\qquad\qquad{}\times \bigl(
\overline{\overline{o}}\bigl(z^2-u^2\bigr) + \overline{
\overline {o}} \bigl((z+u)^2 \bigr) \bigr) \bone_{ \{z\leq0 \}}
\\
&&\quad\qquad{}+ \bigl(\mu(u+z) + uv \sigma^2 - z\zeta\sigma^2 +
\overline {\overline{o}}(v u) + \overline{\overline{o}}(z) \bigr) \bigl( \bigl
\llvert \overline{\overline{o}}\bigl(u^2\bigr) \bigr\rrvert + \bigl
\llvert \overline{\overline{o}}\bigl(z^2\bigr) \bigr\rrvert \bigr)
\bone_{ \{z>0 \}} %
\\
&&\qquad \geq u^2v^2\sigma^4 + 2 uv
\sigma^2 \bigl(\mu(u+z) - z\zeta\sigma^2 + \overline {
\overline{o}}(v u) + \overline{\overline{o}}(z) \bigr) %
\\
&&\quad\qquad{}+ \bigl(\mu(u+z) + uv \sigma^2 - z\zeta\sigma^2 +
\overline {\overline{o}}(v u) + \overline{\overline{o}}(z) \bigr)
\\
&&\qquad\qquad{}\times \bigl(
\overline{\overline{o}}\bigl(z^2-u^2\bigr) + \overline{
\overline {o}} \bigl((z+u)^2 \bigr) \bigr) \bone_{ \{z\leq-u-\varepsilon'' \}}
\\
&&\quad\qquad{}+ v\overline{\overline{o}} \bigl(u^2 \bigr) %
\\
&&\quad\qquad{}+ \bigl(\mu(u+z) + uv \sigma^2 - z\zeta\sigma^2 +
\overline {\overline{o}}(v u) + \overline{\overline{o}}(z) \bigr) \overline{
\overline{o}}\bigl(z^2-u^2\bigr) \bone_{ \{-u+\varepsilon''\leq z\leq0 \}}
\\
&&\quad\qquad{}+ \bigl(uv \sigma^2 + \mu u - z \sigma^2\bigl(\zeta-
\mu/\sigma^2\bigr) + \overline{\overline{o}}(v u) + \overline{
\overline{o}}(z) \bigr)
\\
&&\qquad\qquad{}\times \bigl(\overline{\overline{o}}\bigl(u^2\bigr) +
\overline{\overline{o}}\bigl(z^2\bigr) \bigr) \bone_{ \{z\llvert  \zeta- \mu/\sigma^2\rrvert  > u(v-\varepsilon'') \}},
\end{eqnarray*}
which hold for any (fixed) $\varepsilon''>0$ and all large enough $u$,
$v$ and $\llvert  z\rrvert  $.
In the above, we also make use of the fact that
\begin{eqnarray*}
&&\bigl(\mu(u+z) + uv \sigma^2 - z\zeta\sigma^2 +
\overline {\overline{o}}(v u) + \overline{\overline{o}}(z) \bigr) \bigl(
\overline{\overline{o}}\bigl(z^2-u^2\bigr) + \overline{
\overline {o}} \bigl((z+u)^2 \bigr) \bigr)
\\
&&\qquad{}\times  \bone_{ \{-u-\varepsilon''\leq z\leq-u+\varepsilon'' \}} = v
\overline{\overline{o}} \bigl(u^2 \bigr) %
\end{eqnarray*}
and
\begin{eqnarray*}
&& \bigl( uv \sigma^2 + \mu u - z \sigma^2\bigl(\zeta- \mu/
\sigma^2\bigr) + \overline{\overline{o}}(v u) + \overline{
\overline{o}}(z) \bigr) \bigl( \bigl\llvert \overline{\overline{o}}
\bigl(u^2\bigr) \bigr\rrvert + \bigl\llvert \overline{\overline{o}}
\bigl(z^2\bigr) \bigr\rrvert \bigr)
\\
&&\qquad{}\times  \bone_{ \{0\leq z\llvert  \zeta- \mu/\sigma^2\rrvert  \leq u(v-\varepsilon
'') \}} \geq0
\end{eqnarray*}
hold for all large enough $u$, $v$ and $\llvert  z\rrvert  $.

Finally, choosing $\varepsilon=\varepsilon'=\varepsilon''\in(0,1)$,
we collect the above to obtain
\begin{eqnarray*}
&& \bigl\llvert \psi(v + iu) - \psi(-\zeta-iz) \bigr\rrvert ^2
\\
&&\qquad = \bigl(
\Re \bigl( \psi(v + iu) - \psi(-\zeta-iz) \bigr) \bigr)^2 + \bigl( \Im
\bigl( \psi(v + iu) - \psi(-\zeta-iz) \bigr) \bigr)^2 %
\\
&&\qquad \geq v^4 (1-\varepsilon)\sigma^4/4 +
\bigl(z^2 - u^2 \bigr)^2\sigma^4/4
+ v^2 z^2(1-\varepsilon)\sigma^4/2 -
u^2 v^2 (1+\varepsilon)\sigma^4/2
\\
&&\quad\qquad{}+ u^2v^2\sigma^4 + 2 uv\sigma^2
\bigl(\mu(u+z) - z\zeta\sigma^2 + \overline {\overline{o}}(v u) +
\overline{\overline{o}}(z) \bigr) + v\overline{\overline{o}} \bigl(u^2
\bigr) %
\\
&&\quad\qquad{}+ u^2\overline{\overline{o}} (v ) + \overline{\overline{o}}
\bigl(u^2 \bigr) + u\overline{\overline{o}} \bigl(v^3
\bigr) + v^2\overline{\overline{o}} (u ) + \overline{\overline{o}}
\bigl(v^4 \bigr) %
\\
&&\quad\qquad{}+ \bigl(v^2 + \bigl\llvert z^2 - u^2\bigr
\rrvert \bigr) \bigl( \overline{\overline {o}} \bigl(z^2 -
u^2 \bigr) + \overline{\overline{o}}(uv) + \overline{\overline{o}}
\bigl(v^2 \bigr) + \overline{\overline{o}}(z) \bigr)
\bone_{ \{\llvert  u-\llvert  z\rrvert  \rrvert  >\varepsilon \}} %
\\
&&\quad\qquad{}+ \bigl(\mu(u+z) + uv \sigma^2 - z\zeta\sigma^2 +
\overline {\overline{o}}(v u) + \overline{\overline{o}}(z) \bigr)
\\
&&\qquad\qquad{}\times \bigl(
\overline{\overline{o}}\bigl(z^2-u^2\bigr) + \overline{
\overline {o}} \bigl((z+u)^2 \bigr) \bigr) \bone_{ \{z\leq-u-\varepsilon \}}
\\
&&\quad\qquad{}+ \bigl(\mu(u+z) + uv \sigma^2 - z\zeta\sigma^2 +
\overline {\overline{o}}(v u) + \overline{\overline{o}}(z) \bigr) \overline{
\overline{o}}\bigl(z^2-u^2\bigr) \bone_{ \{-u+\varepsilon\leq z\leq0 \}}
\\
&&\quad\qquad{}+ \bigl(uv \sigma^2 + \mu u - z \sigma^2\bigl(\zeta-
\mu/\sigma^2\bigr) + \overline{\overline{o}}(v u) + \overline{
\overline{o}}(z) \bigr)
\\
&&\qquad\qquad{}\times  \bigl(\overline{\overline{o}}\bigl(u^2\bigr) +
\overline{\overline{o}}\bigl(z^2\bigr) \bigr) \bone_{ \{z\llvert  \zeta- \mu/\sigma^2\rrvert  > u(v-\varepsilon) \}}
\\
&&\qquad \geq v^4 (1-\varepsilon)\sigma^4/4 +
\bigl(z^2 - u^2 \bigr)^2\sigma^4/20
+ v^2 z^2(1-\varepsilon)\sigma^4/2
\\
&&\quad\qquad{}+
u^2 v^2 (1-\varepsilon)\sigma^4/2
+ 2 u^2v\mu\sigma^2 + 2 uvz\sigma^2\bigl(
\mu- \zeta\sigma^2\bigr)
\\
&&\quad\qquad{} + \overline{\overline{o}} \bigl((v
u)^2 \bigr) + uv\overline {\overline{o}}(z) + v\overline{
\overline{o}} \bigl(u^2 \bigr)
\\
&&\quad\qquad{}+ u^2\overline{\overline{o}} (v ) + u\overline{\overline{o}}
\bigl(v^3 \bigr) + v^2\overline{\overline{o}} (u ) +
\overline{\overline{o}} \bigl(v^4 \bigr) %
\\
&&\quad\qquad{}+ \bigl(z^2 - u^2 \bigr)^2\frac{\sigma^4}{20}
\\
&&\quad\qquad{}+ \bigl(v^2 + \bigl\llvert z^2 - u^2\bigr
\rrvert \bigr) \bigl( \overline{\overline{o}} \bigl(z^2 -
u^2 \bigr) + \overline{\overline{o}}(uv) + \overline{\overline{o}}
\bigl(v^2 \bigr) + \overline{\overline{o}}(z) \bigr)
\bone_{ \{\llvert  u-\llvert  z\rrvert  \rrvert  >\varepsilon \}} %
\\
&&\quad\qquad{}+ \bigl(z^2 - u^2 \bigr)^2\frac{\sigma^4}{20}
+ \bigl(\mu(u+z) + uv \sigma^2 - z\zeta\sigma^2 +
\overline{\overline{o}}(v u) + \overline{\overline{o}}(z) \bigr)
\\
&&\hspace*{62pt}\quad\qquad\qquad{}\times  \bigl( \overline{
\overline{o}}\bigl(z^2-u^2\bigr) + \overline{\overline
{o}} \bigl((z+u)^2 \bigr) \bigr) \bone_{ \{z\leq-u-\varepsilon \}} %
\\
&&\quad\qquad{}+ \bigl(z^2 - u^2 \bigr)^2\frac{\sigma^4}{20}
\\
&&\quad\qquad{}
+ \bigl(\mu(u+z) + uv \sigma^2 - z\zeta\sigma^2 +
\overline{\overline{o}}(v u) + \overline{\overline{o}}(z) \bigr) \overline{
\overline{o}}\bigl(z^2-u^2\bigr) \bone_{ \{-u+\varepsilon\leq z\leq0 \}}
\\
&&\quad\qquad{}+ \bigl(z^2 - u^2 \bigr)^2\frac{\sigma^4}{20}
+ \bigl(uv \sigma ^2 + \mu u - z \sigma^2\bigl(\zeta-
\mu/\sigma^2\bigr) + \overline{\overline {o}}(v u) + \overline{
\overline{o}}(z) \bigr)
\\
&&\hspace*{117pt}{}\times  \bigl(\overline{\overline{o}}\bigl(u^2\bigr) +
\overline{\overline{o}}\bigl(z^2\bigr) \bigr) \bone_{ \{z\llvert  \zeta- \mu/\sigma^2\rrvert  > u(v-\varepsilon) \}}.
\end{eqnarray*}
Let us estimate the above terms separately:
%
\begin{eqnarray}
&& v^4 (1-\varepsilon)\sigma^4/4 + \bigl(z^2
- u^2 \bigr)^2\sigma^4/20 + v^2
z^2(1-\varepsilon)\sigma^4/2 + u^2
v^2 (1-\varepsilon)\sigma^4/2 \nonumber
\\
&&\quad{}+ 2 u^2v\mu\sigma^2 + 2 uvz\sigma^2\bigl(
\mu- \zeta\sigma^2\bigr) + \overline{\overline{o}} \bigl((v
u)^2 \bigr) + uv\overline {\overline{o}}(z) + v\overline{
\overline{o}} \bigl(u^2 \bigr) \nonumber
\\
&&\quad{}+ u^2\overline{\overline{o}} (v ) + u\overline{\overline{o}}
\bigl(v^3 \bigr) + v^2\overline{\overline{o}} (u ) +
\overline{\overline{o}} \bigl(v^4 \bigr) \nonumber
\\
&&\qquad \geq v^4 (1-\varepsilon-\delta)\sigma^4/4 +
\bigl(z^2 - u^2 \bigr)^2\sigma^4/20
+ v^2 z^2(1-\varepsilon)\sigma^4/2 \nonumber
\\
&&\quad\qquad{} +
u^2 v^2 (1-\varepsilon- \delta)\sigma^4/2+ u\overline{\overline{o}} \bigl(v^3 \bigr) - c_1 uv
\llvert z\rrvert,\nonumber
\\
&& \bigl(z^2 - u^2 \bigr)^2
\frac{\sigma^4}{20} + \bigl(v^2 + \bigl\llvert z^2 -
u^2\bigr\rrvert \bigr)\nonumber
\\
 &&\hspace*{71pt}{}\times \bigl( \overline{\overline{o}}
\bigl(z^2 - u^2 \bigr) + \overline{\overline{o}}(uv) +
\overline{\overline{o}} \bigl(v^2 \bigr) + \overline{\overline{o}}(z)
\bigr) \bone_{ \{\llvert  u-\llvert  z\rrvert  \rrvert  >\varepsilon \}}\nonumber
\\
\label{eqaux1}&&\qquad \geq \bigl(z^2 - u^2 \bigr)^2 \biggl(
\frac{\sigma^4}{20} + \biggl(\frac{v^2}{\llvert  z\rrvert   + u} + 1 \biggr) \frac{\overline{\overline{o}} (\llvert  z\rrvert   + u  ) + \overline
{\overline{o}}(uv)
+ \overline{\overline{o}} (v^2 ) + \overline{\overline
{o}}(z)}{\llvert  z\rrvert  +u}
\biggr),
\\
&&  \bigl(z^2 - u^2 \bigr)^2
\frac{\sigma^4}{20} + \bigl(\mu(u+z) + uv \sigma^2 - z\zeta
\sigma^2 + \overline{\overline{o}}(v u) + \overline{\overline{o}}(z)
\bigr)\nonumber
\\
&&\hspace*{71pt}{}\times  \bigl( \overline{\overline{o}}\bigl(z^2-u^2\bigr) +
\overline{\overline {o}} \bigl((z+u)^2 \bigr) \bigr)
\bone_{ \{z\leq-u-\varepsilon \}}\nonumber
\\
\label{eqaux2}&&\qquad \geq \bigl(z^2 - u^2 \bigr)^2 \biggl(
\frac{\sigma^4}{20} + \frac{\mu(u+z) + uv \sigma^2 - z\zeta\sigma^2 + \overline
{\overline{o}}(v u) + \overline{\overline{o}}(z)}{u+\llvert  z\rrvert  }
\\
&&\hspace*{166pt}{}\times \frac{\overline{\overline{o}}(z^2-u^2) + \overline{\overline
{o}} ((\llvert  z\rrvert  -u)^2 )}{(\llvert  z\rrvert  -u)\llvert  z^2-u^2\rrvert  } \biggr), \nonumber
\\
&& \bigl(z^2 - u^2 \bigr)^2
\frac{\sigma^4}{20} + \bigl(\mu(u+z) + uv \sigma^2 - z\zeta
\sigma^2 + \overline{\overline{o}}(v u) + \overline{\overline{o}}(z)
\bigr) \nonumber
\\
&&\hspace*{70pt}{}\times \overline{\overline{o}}\bigl(z^2-u^2\bigr)
\bone_{ \{-u+\varepsilon\leq z\leq0 \}}\nonumber
\\
\label{eqaux3}
&&\qquad \geq \bigl(z^2 - u^2 \bigr)^2
\\
&&\quad\qquad{}\times  \biggl(
\frac{\sigma^4}{20} + \frac{ (\mu(u+z) + uv \sigma^2 - z\zeta\sigma^2 + \overline
{\overline{o}}(v u) + \overline{\overline{o}}(z) )}{\llvert  z\rrvert   + u} \frac{\overline{\overline{o}}(\llvert  z\rrvert   + u)}{(u-\llvert  z\rrvert  )^2} \biggr), \nonumber
\\
&& \bigl(z^2 - u^2 \bigr)^2
\frac{\sigma^4}{20} + \bigl(uv \sigma^2 + \mu u - z
\sigma^2\bigl(\zeta- \mu/\sigma^2\bigr) + \overline{
\overline {o}}(v u) + \overline{\overline{o}}(z) \bigr)\nonumber
\\
&&\hspace*{71pt}{}\times  \bigl(\overline{
\overline{o}}\bigl(u^2\bigr) + \overline{\overline{o}}
\bigl(z^2\bigr) \bigr) \bone_{ \{z\llvert  \zeta- \mu/\sigma^2\rrvert  > u(v-\varepsilon) \}}\nonumber
\\
\label{eqaux4}
&&\qquad \geq \bigl(z^2 - u^2 \bigr)^2
\\
&&\qquad\quad{}\times \biggl(
\frac{\sigma^4}{20} + \frac{uv \sigma^2 + \mu u - z \sigma^2(\zeta- \mu/\sigma^2) +
\overline{\overline{o}}(v u) + \overline{\overline{o}}(z)}{\llvert  z\rrvert  +u} \frac{\overline{\overline{o}}(z^2)}{z^2} \biggr),\nonumber
\end{eqnarray}
where we fix arbitrary $\delta\in(0,1-\varepsilon)$ and assume that
$u$, $v$ and $\llvert  z\rrvert  $ are large enough.

It only remans to notice that for any $v>0$, there exist $R_1>0$ and
$R_2>0$, such that for all $\llvert  z\rrvert  >R_1$ and $u>R_2$, the right-hand sides
of (\ref{eqaux1})--(\ref{eqaux4}) are nonnegative, and in addition,
\[
u^2 v^2 (1-\varepsilon- \delta)\sigma^4/2 +
u\overline{\overline{o}} \bigl(v^3 \bigr) \geq u^2
v^2 (1-\varepsilon- \delta)\sigma^4/4. %
\]
Thus we conclude that for any large enough $v>0$, there exist $R_1>0$,
$R_2>0$ and $ \{ c_i>0 \}$, such that, for all $\llvert  z\rrvert  >R_1$ and
$u>R_2$, the following holds:
\begin{eqnarray*}
&& \bigl\llvert \psi(v + iu) - \psi(-\zeta-iz) \bigr\rrvert ^2
\\
&&\qquad \geq
c_{2} v^4 + c_{3}\bigl(z^2 -
u^2 \bigr)^2 + 2c_{4} v^2
z^2 + c_{5} u^2 v^2 -
c_{1} u v \llvert z\rrvert %
\\
&&\qquad = c_{2} v^4 + c_{3}\bigl(z^2 -
u^2 \bigr)^2 + c_{4} v^2
z^2 + \bigl(\sqrt{c_{4}} v \llvert z\rrvert -
\sqrt{c_{5}} u v\bigr)^2
\\
&&\quad\qquad{} + u v \llvert z\rrvert (2
\sqrt{c_4 c_5} v - c_{1}) %
\\
&&\qquad\geq c_{6} \bigl( \bigl(z^2 - u^2
\bigr)^2 + z^2 \bigr). %
\end{eqnarray*}\upqed
\end{pf*}\vspace*{-9pt}

\section{}\label{appB}\vspace*{-9pt}
\begin{pf*}{Proof of Theorem~\ref{thmain}}
The proof consists of four steps. In step~1, using only Assumptions
\ref{ass1} and~\ref{ass3}, we construct the generalized function
$g\in\mathcal{D}^*$ as a limit of $g_r$, to show that it has support
in $[0,\infty)$ and to establish the rate of convergence~(\ref
{eqmainrategrweakpt}). In step~2, assuming in addition that $\hat
{h}\in\mathbb{L}^1(\RR)$, we show that $g$ coincides with a
continuous function in $(0,\infty)$, and that $g_r$ converges
pointwise, with the rate of convergence given in the theorem.
In step~3, we make the additional assumption~(\ref{eqtailcond}) to prove
that $g$ is locally integrable and that $\EE g(X_t) = \EE h(X_t)$ for
all $t>0$. Finally, in step~4, we use the results of steps 2~and~3, to
show that, even in\vadjust{\goodbreak} the absence of additional assumptions on $\hat{h}$
(i.e., using only Assumptions~\ref{ass1} and~\ref{ass3}), the
generalized function $g\in\mathcal{D}^*$, constructed in step~1,
satisfies $\langle g, p_t \rangle= \EE h(X_t)$, for all $t>0$.

\begin{longlist}[\textit{Step} 2.]
\item[\textit{Step} 1.] First, for any $w\geq\gamma$, we introduce
%
\begin{equation}
\label{eqgwrdef} \qquad g^w_r(x) = \frac{1}{2\pi i} \int
_{\mathcal{G}^w_r} e^{\lambda x} \int_{\RR} \biggl(
\frac{\psi'(\lambda)}{\psi(\lambda) - \psi
(-\zeta-iz)} - \frac{1}{\lambda+ \zeta+ iz} \biggr) \hat{h}(z) \,dz \,d\lambda,
\end{equation}
%
where we denote by $\mathcal{G}^w_r$ the vertical interval $[w-ir,w+ir]$.
Notice that $g_r=g^{\gamma}_r$.
Let us show that $g^w_r$ has a weak limit $g\in\mathcal{D}^*$, as
$r\rightarrow\infty$, which is independent of~$w$.
Fix any $r'>r$, and any test function $\phi\in\mathcal{D}$, and
proceed as follows:
%
\begin{eqnarray}\label{eqgwract}
&& \bigl\llvert \bigl\langle g^{w}_{r'} -
g^w_{r},\phi\bigr\rangle\bigr\rrvert\nonumber
\\
&&\qquad  \leq
\frac{1}{2\pi} \int_{r<\llvert  u\rrvert  <r'} \biggl\llvert \int
_{\RR} e^{(w + iu) x} \phi(x) \,dx\biggr\rrvert
\nonumber\\[-8pt]\\[-8pt]\nonumber
&&\hspace*{86pt}{}\times  \biggl\llvert
\int_{\RR} \biggl( \frac{\psi'(w + iu)}{\psi(w + iu) - \psi
(-\zeta-iz)}
\\
&&\hspace*{150pt}{} - \frac{1}{w + iu + \zeta+ iz}
\biggr) \hat{h}(z) \,dz\biggr\rrvert \,du. \nonumber
\end{eqnarray}
Integrating by parts repeatedly, we obtain
%
\begin{equation}
\label{eqephiintest} \biggl\llvert \int_{\RR} e^{(w + iu) x}
\phi(x) \,dx\biggr\rrvert = \llvert w + iu \rrvert ^{-3} \int
_{\RR} e^{w x} \bigl\llvert \phi '''(x)
\bigr\rrvert \,dx.
\end{equation}
%
Making use of (\ref{eqcor2est1}), we obtain
\[
\bigl\llvert \bigl\langle g^w_{r'}- g^w_{r},
\phi\bigr\rangle\bigr\rrvert \leq \frac{c_1}{\pi} \bigl( 1/r -
1/r' \bigr) \bigl(\llVert \hat{h}\rrVert _{\mathbb{L}^1(\RR)}\wedge
\llVert \hat{h}\rrVert _{\mathbb
{L}^2(\RR)} \bigr) \int_{\RR}
e^{w x} \bigl\llvert \phi'''(x)
\bigr\rrvert \,dx. %
\]
The above estimate shows that $g^w = \lim_{r\rightarrow\infty}g^w_r$
is well defined as an element of~$\mathcal{D}^*$. To see that $g^w$ is
independent of $w$, consider arbitrary $w'>w$ and connect the two
intervals of integration $\mathcal{G}^w_r$ and $\mathcal{G}^{w'}_r$
by the two horizontal parts: $\mathcal{C}^1_r = [w-ir,w'-ir]$ and
$\mathcal{C}^2_r = [w+ir,w'+ir]$. Since the integrand in (\ref
{eqgwrdef}) is analytic with respect to $\lambda$, the integral over
the closed contour (with appropriately chosen directions on each part)
is zero. Thus we only need to show that the integrals over $\mathcal
{C}^1_r$ and $\mathcal{C}^2_r$ vanish, as $r\rightarrow\infty$,
\begin{eqnarray*}
&& \bigl\llvert \bigl\langle g^{w'}_{r} - g^w_{r},
\phi\bigr\rangle\bigr\rrvert
\\
&&\qquad \leq
\frac{1}{2\pi} \sum_{u=-r,r} \int
_{w}^{w'} \biggl\llvert \int_{\RR}
e^{(v + iu) x} \phi(x) \,dx\biggr\rrvert
\\
&&\hspace*{97pt}{}\times \biggl\llvert \int_{\RR}
\biggl( \frac{\psi'(v + iu)}{\psi(v + iu) - \psi
(-\zeta-iz)}
\\
&&\hspace*{161pt}{} - \frac{1}{v + iu + \zeta+ iz} \biggr) \hat{h}(z) \,dz\biggr
\rrvert \,dv. %
\end{eqnarray*}
Estimates (\ref{eqephiintest}) and (\ref{eqcor2est1}) imply
that the right-hand side of the above vanishes as $r\rightarrow\infty
$. Thus $g^w$ is independent of $w$, and we denote it by $g$.
Let us show that $g$ has support in $[0,\infty)$. Choose an arbitrary
$\phi\in\mathcal{D}$, such that $\operatorname{supp} (\phi)
\subset(-\infty,0)$, and consider $\langle g^{w}_{r},\phi\rangle$.
Equation (\ref{eqephiintest}), in this case, becomes
\[
\biggl\llvert \int_{-\infty}^0 e^{(w + iu) y}
\phi(y) \,dy\biggr\rrvert = \llvert w + iu \rrvert ^{-3} \int
_{-\infty}^0 e^{w y} \bigl\llvert \phi
'''(y)\bigr\rrvert \,dy \leq
c_2 \llvert w + iu \rrvert ^{-3}, %
\]
which holds uniformly over all $u\in\RR$ and $w\geq\gamma$.
Thus we can close the contour of integration in the integral
representation of $\langle g^{w}_{r},\phi\rangle$ [cf. (\ref
{eqgwract})] by a semicircle (on the right-hand side), and using the
above estimate, along with (\ref{eqcor2est1}) and the analyticity of
the integrand in $H_R$ (cf. Corollary~\ref{coranExt}), conclude that
$\langle g^{w}_{r},\phi\rangle\rightarrow0$, as $r\rightarrow\infty$.
To obtain (\ref{eqmainrategrweakpt}), we recall (\ref
{eqgwract}) and the fact that
\[
\int_{\RR} e^{(w + iu) x} p_t(x) \,dx = \exp
\bigl(t\psi(w+iu) \bigr). %
\]
Then (\ref{eqmainrategrweakpt}) follows from (\ref
{eqcor2est1}) and (\ref{eqfourdensasymp}).


\item[\textit{Step} 2.] Next, under the additional assumption that $\hat{h}$ is
absolutely integrable, we show that $g$ coincides with a continuous
function in $(0,\infty)$, and that $g_r(x)$ converges to $g(x)$ for
every $x>0$.
Applying Fubini's theorem and integration by parts we obtain
%
\begin{eqnarray}\label{eqLevyfinalgRconv}
&& 2\pi i g_r(x) \nonumber
\\
&&\qquad = \frac{1}{2\pi i} \int_{\mathcal{G}_r}
e^{\lambda x} \int_{\RR} \biggl( \frac{\psi'(\lambda)}{\psi(\lambda) - \psi
(-\zeta-iz)} -
\frac{1}{\lambda+ \zeta+ iz} \biggr) \hat{h}(z) \,dz \,d\lambda \nonumber
\\
&&\qquad = \int_{\mathcal{G}_r} \int_{\RR} e^{\lambda x}
\frac{\psi '(\lambda)}{\psi(\lambda) - \psi(-\zeta-iz)} \hat{h}(z) \,dz \,d\lambda \nonumber
\\
&&\quad\qquad{}- \int_{\mathcal{G}_r} \int
_{\RR} e^{\lambda x} \frac{1}{\lambda+
\zeta+ iz} \hat{h}(z) \,dz \,d
\lambda \nonumber
\\
&&\qquad = \frac{1}{x} \int_{\RR} \int_{\mathcal{G}_r}
\bigl(e^{\lambda
x} \bigr)' \frac{\psi'(\lambda)}{\psi(\lambda) - \psi(-\zeta
-iz)} \,d\lambda
\hat{h}(z) \,dz \nonumber
\\
&&\quad\qquad{}- \frac{1}{x} \int_{\RR} \int
_{\mathcal{G}_r} \bigl( e^{\lambda
x} \bigr)'
\frac{1}{\lambda+ \zeta+ iz} \,d\lambda\hat{h}(z) \,dz \nonumber
\\
&&\qquad = \frac{1}{x} \int_{\RR} \biggl(
\frac{\exp((\gamma+ ir)x) \psi
'(\gamma+ ir)}{\psi(\gamma+ ir) - \psi(-\zeta-iz)} - \frac{\exp((\gamma- ir)x) \psi'(\gamma- ir)}{\psi(\gamma- ir)
- \psi(-\zeta-iz)} \biggr) \hat{h}(z) \,dz\nonumber
\\
&&\qquad\quad{} + \frac{1}{x} \int_{\RR} \biggl(\frac{\exp((\gamma- ir)x)}{\gamma
+ \zeta- ir + iz} -
\frac{\exp((\gamma+ ir)x)}{\gamma+ \zeta+ ir + iz} \biggr) \hat {h}(z) \,dz \nonumber
\\
&&\quad\qquad{}- \frac{1}{x} \int_{\RR} \int_{\mathcal{G}_r}
\frac{\exp(\lambda
x)\psi''(\lambda)}{\psi(\lambda) - \psi(-\zeta-iz)} \,d\lambda\hat {h}(z) \,dz\nonumber
\\
&&\quad\qquad{} + \frac{1}{x} \int
_{\RR} \int_{\mathcal{G}_r} \frac{\exp(\lambda
x) (\psi'(\lambda) )^2}{ (\psi(\lambda) - \psi
(-\zeta-iz) )^2} \,d
\lambda\hat{h}(z) \,dz
\\
&&\quad\qquad{}- \frac{1}{x} \int_{\RR} \int_{\mathcal{G}_r}
\frac{\exp(\lambda
x)}{(\lambda+ \zeta+ iz)^2} \,d\lambda\hat{h}(z) \,dz. \nonumber
\end{eqnarray}
Let us show that the first integral on the right-hand side of (\ref
{eqLevyfinalgRconv}) converges to zero, as $r\rightarrow\infty$.
Due to (\ref{eqpsipsiprimeasymp}), for any $R>0$, there exist $c_1,
c_2, c_3>0$ and $r'>0$, such that, for all $r>r'$, we have
\begin{eqnarray}\label{eqthaux1}
&& \biggl\llvert \int_{\RR} \biggl(
\frac{\exp((\gamma+ ir)x) \psi'(\gamma
+ ir)}{\psi(\gamma+ ir) - \psi(-\zeta-iz)} - \frac{\exp((\gamma- ir)x) \psi'(\gamma- ir)}{\psi(\gamma- ir)
- \psi(-\zeta-iz)} \biggr) \hat{h}(z) \,dz \biggr\rrvert \nonumber
\\
&&\qquad \leq
c_1 e^{\gamma x} \int_{\llvert  z\rrvert  \leq R}
\frac{r}{r^2 - c_2} \bigl\llvert \hat {h}(z)\bigr\rrvert \,dz
\nonumber\\[-8pt]\\[-8pt]\nonumber
&&\quad\qquad{}
+ c_3 e^{\gamma x}\int_{\llvert  z\rrvert  >R} \biggl(
\frac{ r}{\llvert  \psi(\gamma+
ir) - \psi(-\zeta-iz)\rrvert  }\nonumber
\\
&&\hspace*{106pt}{} + \frac{r}{\llvert  \psi(\gamma- ir) - \psi(-\zeta-iz) \rrvert  } \biggr) \bigl\llvert \hat{h}(z) \bigr
\rrvert \,dz. \nonumber
\end{eqnarray}
We choose $R$ to be large enough, so that the estimate in Lemma~\ref
{lepsiest2} can be applied for all $\llvert  z\rrvert  \geq R$ and $\llvert  u\rrvert  =r$,
%
\begin{eqnarray*}
&& \frac{1}{\llvert  \psi(\gamma+ ir) - \psi(-\zeta-iz) \rrvert  } + \frac{1}{\llvert  \psi(\gamma- ir) - \psi(-\zeta-iz) \rrvert  }
\\
&&\qquad \leq c_4
\frac{1}{\sqrt{(z^2 - r^2)^2 + z^2}}. 
\end{eqnarray*}
Thus, the first integral on the right-hand side of (\ref
{eqLevyfinalgRconv}) is estimated from above by
%
\begin{eqnarray}\label{eqfirstintrate}
&& c_1 e^{\gamma x} \int_{\llvert  z\rrvert  \leq N}
\frac{r}{r^2 - c_2} \bigl\llvert \hat{h}(z) \bigr\rrvert \,dz + c_3
c_4 e^{\gamma x} \int_{\llvert  z\rrvert  >N}
\frac{r}{\sqrt{(z^2 - r^2)^2 +
z^2}} \bigl\llvert \hat{h}(z) \bigr\rrvert \,dz\nonumber
\\
&&\qquad \leq c_1 e^{\gamma x} \int_{\llvert  z\rrvert  \leq N}
\frac{r}{r^2 - c_2} \bigl\llvert \hat {h}(z)\bigr\rrvert \,dz + c_3c_4
e^{\gamma x} \int_{\llvert  z\rrvert  >N,\llvert  z\rrvert  \leq r/2} \frac{2}{r} \bigl\llvert
\hat {h}(z)\bigr\rrvert \,dz
\\
&&\qquad\quad{} + 2c_3c_4 e^{\gamma x}
\int_{\llvert  z\rrvert  >N,\llvert  z\rrvert  >r/2} \bigl\llvert \hat{h}(z) \bigr\rrvert \,dz,\nonumber
\end{eqnarray}
which vanishes, as $r\rightarrow\infty$.
Similarly, we proceed with the second integral on the right-hand side
of (\ref{eqLevyfinalgRconv}):
%
\begin{eqnarray}\label{eqsecondintrate}
&& \biggl\llvert \int_{\RR} \biggl(\frac{\exp((\gamma- ir)x)}{\gamma+ \zeta
- ir + iz} -
\frac{\exp((\gamma+ ir)x)}{\gamma+ \zeta+ ir + iz} \biggr) \hat {h}(z) \,dz \biggr\rrvert \nonumber
\\
&&\qquad \leq2 e^{\gamma x}
\int_{\RR} \frac{1}{\sqrt{(\gamma+ \zeta)^2 +
(z - r)^2}} \bigl\llvert \hat{h}(z) \bigr
\rrvert \,dz %
\\
&&\qquad \leq2 e^{\gamma x} \int_{\llvert  z\rrvert  \leq r/2}
\frac{2}{r} \bigl\llvert \hat{h}(z) \bigr\rrvert \,dz + 2 e^{\gamma x}
\int_{\llvert  z\rrvert  > r/2} \frac{1}{\gamma+ \zeta} \bigl\llvert \hat {h}(z)\bigr
\rrvert \,dz.\nonumber
\end{eqnarray}
As $r\rightarrow\infty$, the third, fourth and fifth integrals on the
right-hand side of (\ref{eqLevyfinalgRconv}) converge uniformly
over $x$, changing on any compact in $(0,\infty)$.
Let us prove the convergence of the third integral.
We consider arbitrary, large enough, $r^{\prime}>r$ and proceed as follows:
%
\begin{eqnarray}
\label{eqLevyfinalgRconv2}
&& \biggl\llvert \int_{\RR} \int_{\mathcal{G}_{r'}\setminus\mathcal{G}_r}
\frac{\exp(\lambda x)\psi''(\lambda)}{\psi(\lambda) - \psi
(-\zeta-iz)} \,d\lambda\hat{h}(z) \,dz\biggr\rrvert \nonumber
\\
&&\qquad \leq2e^{\gamma x} \int
_{\RR} \int_{\llvert  u\rrvert  \in[r,r']} \frac{\llvert  \psi
''(\gamma+ iu)\rrvert  }{\llvert  \psi(\gamma+ iu) - \psi(-\zeta-iz)\rrvert  } \,du
\bigl\llvert \hat {h}(z)\bigr\rrvert \,dz
\nonumber\\[-8pt]\\[-8pt]\nonumber
&&\qquad \leq c_5 e^{\gamma x} \int_{-R}^R
\bigl\llvert \hat{h}(z) \bigr\rrvert \,dz \int_{r}^{r^{\prime}}
\frac{1}{u^2 - c_6} \,du
\\
&&\quad\qquad{}+ c_5 e^{\gamma x} \int
_{\llvert  z\rrvert  >R} \int_{r}^{r'} \biggl
\llvert \frac{1}{\psi(\gamma+ iu) - \psi(-\zeta-iz)} \biggr\rrvert \,du \bigl\llvert \hat{h}(z) \bigr\rrvert
\,dz,\nonumber
\end{eqnarray}
where, again, we choose $R>0$ to be large enough, so that Lemma~\ref
{lepsiest2} can be applied for all $\llvert  z\rrvert  \geq R$ and all $\llvert  u\rrvert  \geq r$.
It is easy to see that the first term on the right-hand side of (\ref
{eqLevyfinalgRconv2}) vanishes, as $r,r'\rightarrow\infty$.
To show that the second term on the right-hand side of (\ref
{eqLevyfinalgRconv2}) vanishes, as $r,r'\rightarrow\infty$,
we make use of Lemma~\ref{lepsiest2} to obtain
\begin{eqnarray*}
&& \int_{\llvert  z\rrvert  >R} \int_{r}^{r'}
\biggl\llvert \frac{1}{\psi(\gamma+ iu) - \psi(-\zeta-iz)} \biggr\rrvert \,du \bigl\llvert \hat{h}(z) \bigr
\rrvert \,dz
\\
&&\qquad \leq c_7 \int_{\llvert  z\rrvert  >R} \int
_{r}^{r'} \frac{1}{\sqrt{(u^2-z^2)^2 + z^2}} \,du \bigl\llvert
\hat{h}(z) \bigr\rrvert \,dz. %
\end{eqnarray*}
Thus the right-hand side of (\ref{eqLevyfinalgRconv2}) is bounded
from above by the following expression:
%
\begin{eqnarray}
\label{eqthirdintrate}
&& c_8 e^{\gamma x} \frac{1}{r} \int_{-R}^R
\bigl\llvert \hat{h}(z) \bigr\rrvert \,dz + c_9 e^{\gamma x} \int
_{\llvert  z\rrvert  >R} \int_{r}^{\infty}
\frac{1}{\sqrt{(u^2-z^2)^2 + z^2}} \,du \bigl\llvert \hat{h}(z) \bigr\rrvert \,dz \nonumber
\\
&&\qquad \leq c_8 e^{\gamma x} \frac{1}{r} \int
_{-R}^R \bigl\llvert \hat{h}(z) \bigr\rrvert \,dz
\nonumber\\[-8pt]\\[-8pt]\nonumber
&&\quad\qquad{} + c_9 e^{\gamma x} \int_{R<\llvert  z\rrvert  \leq r/2} \int
_{r}^{\infty} \frac{1}{(u-r/2)^2 } \,du \bigl\llvert
\hat{h}(z) \bigr\rrvert \,dz
\\
&&\quad\qquad{}
+ c_9 e^{\gamma x} \int_{\llvert  z\rrvert  >R,\llvert  z\rrvert  > r/2} \int
_{0}^{\infty} \frac{1}{\sqrt{u^4+1}} \,du \bigl\llvert
\hat{h}(z) \bigr\rrvert \,dz, \nonumber
\end{eqnarray}
where we make use of
%
\begin{equation}
\label{eqthirdintaux} \biggl\llvert \int_{r}^{r'}
\frac{1}{\sqrt{(u^2-z^2)^2 + z^2}} \,du\biggr\rrvert \leq\int_{\RR}
\frac{1}{\sqrt{u^4 + 1}} \,du < \infty.
\end{equation}
Next, we use Lemma~\ref{lepsiest2} to prove the convergence of
the fourth integral on the right-hand side of (\ref{eqLevyfinalgRconv}):
%
\begin{eqnarray}\label{eqfourthint}
&& \int_{\RR} \int_{r}^{r^{\prime}}
\frac{\llvert  \exp((\gamma+iu)x)\rrvert  \llvert  \psi'(\gamma+ iu)\rrvert  ^2}{\llvert   \psi(\gamma+ iu) - \psi(-\zeta-iz)\rrvert  ^2} \,du \bigl\llvert \hat {h}(z)\bigr\rrvert \,dz\nonumber
\\
&&\qquad
\leq c_{10} e^{\gamma x} \int_{\llvert  z\rrvert  \leq R} \int
_{r}^{r^{\prime}} \frac{u^2}{\llvert   u^2 - c_{11}\rrvert  ^2} \,du \bigl\llvert
\hat{h}(z) \bigr\rrvert \,dz
\\
&&\quad\qquad{}  + c_{12}e^{\gamma x} \int
_{\llvert  z\rrvert  >R} \int_{r}^{r'}
\frac{u^2}{(u^2 - z^2
)^2 + z^2} \,du \bigl\llvert \hat{h}(z) \bigr\rrvert \,dz. \nonumber
\end{eqnarray}
The right-hand side of the above is bounded by
%
\begin{eqnarray}\label{eqfourthintrate}
\qquad&& c_{13} e^{\gamma x}\frac{1}{r} \int_{\llvert  z\rrvert  \leq R}
\bigl\llvert \hat{h}(z) \bigr\rrvert \,dz
+ c_{12}e^{\gamma x} \int
_{R<\llvert  z\rrvert  \leq r/2} \int_{r}^{\infty}
\frac
{u^2}{(u - r/2 )^4} \,du \bigl\llvert \hat{h}(z) \bigr\rrvert \,dz %
\nonumber\\[-8pt]\\[-8pt]\nonumber
&&\qquad{} + c_{12}e^{\gamma x} \int_{R\langle \llvert  z\rrvert  \rangle r/2}
\biggl( 2 (1+1/R )^2 + \int_{\llvert  u\rrvert  >1}
\frac{1}{u^2} \,du \biggr) \bigl\llvert \hat{h}(z) \bigr\rrvert \,dz,\nonumber
\end{eqnarray}
where we make use of
%
\begin{eqnarray}\label{eqfourthintaux}
&& \biggl\llvert \int_{0}^{\infty} \frac{u^2}{(u^2 - z^2 )^2 + z^2}
\,du \biggr\rrvert\nonumber
\\
&&\qquad  \leq\int_{z-1}^{z+1}
\frac{u^2}{(u^2 - z^2)^2 + z^2} \,du
+\int_{u\in[0,z-1]\cup[z+1,\infty)} \frac{u^2}{(u^2 - z^2)^2 +z^2} \,du
\nonumber
\\
&&\qquad \leq\int_{1-1/z}^{1+1/z}
\frac{zu^2}{z^2(u^2 - 1)^2 + 1} \,du
\\
&&\quad\qquad{} + \int_{u\in[0,z-1]\cup[z+1,\infty)} \frac{u^2}{(u - z)^2 (u+z)^2} \,du\nonumber
\\
&&\qquad \leq2 (1+1/R )^2 + \int_{\llvert  u\rrvert  >1}
\frac{1}{u^2} \,du < \infty.\nonumber
\end{eqnarray}
To show that the last integral on the right-hand side of (\ref
{eqLevyfinalgRconv}) converges as $r\rightarrow\infty$, we notice that
%
\begin{eqnarray}\label{eqfifthint}
&& \biggl\llvert \int_{\RR} \int
_{\mathcal{G}_{r'}\setminus\mathcal{G}_r} \frac{\exp(\lambda x)}{(\lambda+ \zeta+ iz)^2} \,d\lambda\hat{h}(z) \,dz \biggr\rrvert\nonumber
\\
&&\qquad
\leq e^{\gamma x} \int_{\RR} \int_{r}^{r'}
\frac{1}{(\gamma+ \zeta
)^2 + (z+u)^2} \,du \bigl\llvert \hat{h}(z) \bigr\rrvert \,dz
\nonumber\\[-8pt]\\[-8pt]\nonumber
&&\qquad \leq e^{\gamma x} \int_{\llvert  z\rrvert  \leq r/2} \int_{r}^{\infty}
\frac{1}{(u-r/2)^2} \,du \bigl\llvert \hat{h}(z) \bigr\rrvert \,dz
\\
&&\quad\qquad{}  + e^{\gamma x}
\int_{\llvert  z\rrvert  > r/2} \int_{\RR} \frac{1}{(\gamma+ \zeta
)^2 + u^2}
\,du \bigl\llvert \hat{h}(z) \bigr\rrvert \,dz. \nonumber
\end{eqnarray}
Thus we have shown that for $x>0$, $g_r(x)$ converges to $g(x)$, as
$r\rightarrow\infty$. Moreover, estimates (\ref{eqfirstintrate}),
(\ref{eqsecondintrate}), (\ref{eqthirdintrate}), (\ref
{eqfourthintrate}) and (\ref{eqfifthint}) imply the desired rate
of convergence (\ref{eqmainrategr}).
To see that the limiting function $g(x)$ is continuous for $x\in
(0,\infty)$, we notice that, due to (\ref{eqmainrategr}), the
convergence is uniform over $x$ changing on any compact in $(0,\infty
)$. Similarly, (\ref{eqmainrategr}) implies that $g$ has at most
exponential growth at infinity.

\item[\textit{Step} 3.] Under the additional assumption that $\hat{h}$ is
absolutely integrable and that~(\ref{eqtailcond}) holds, let us show
that for every $w\geq\gamma$, the function
\[
F^w\dvtx  u\mapsto\int_{\RR} \biggl(
\frac{\psi'(w+iu)}{\psi(w+iu) -
\psi(-\zeta-iz)} - \frac{1}{wa+iu + \zeta+ iz} \biggr) \hat{h}(z) \,dz %
\]
is square integrable over $\RR$.
Notice that $F^w$ is continuous (cf. Corollary~\ref{coranExt}), and
therefore, it suffices to estimate $\llvert  F^w(u)\rrvert  $ for large $\llvert  u\rrvert  $. Let us
choose a large enough $R>0$, for which there exists a constant $\beta
\in(0,1)$, such that
\[
\bigl\llvert \psi(w+iu)\bigr\rrvert \geq2 \bigl\llvert \psi(-\zeta-iz)\bigr
\rrvert %
\]
holds for all $\llvert  z\rrvert  \leq\beta\llvert  u\rrvert  $ and all $\llvert  u\rrvert  >R$. Notice that such $R$
and $\beta$ do exist, due to the inequalities (\ref
{eqpsipsiprimeasymp}). Increasing $R$, if necessary, we ensure that
Lemma~\ref{lepsiest2} can be applied for all $\llvert  u\rrvert  >R$ and all
$\llvert  z\rrvert  >\beta R$.
Finally, for any $\llvert  u\rrvert  >R$, we integrate by parts to obtain
\begin{eqnarray*}
\hspace*{-1pt}\bigl\llvert F^w(u)\bigr\rrvert &\leq&\int_{\llvert  z\rrvert  \leq\beta\llvert  u\rrvert  }
\biggl\llvert \frac{\psi'(w+iu)}{\psi(w+iu)
- \psi(-\zeta-iz)} - \frac{1}{w+iu + \zeta+ iz} \biggr\rrvert \bigl\llvert
\hat{h}(z)\bigr\rrvert \,dz %
\\[-1pt]
\hspace*{-1pt}&&{}+\sum_{z=-\beta u,\beta u} \biggl\llvert \frac{\psi'(w+iu)}{\psi(w+iu) -
\psi(-\zeta-iz)}
\\[-1pt]
\hspace*{-1pt}&&\hspace*{85pt}{}-
\frac{1}{w+iu + \zeta+ iz } \biggr\rrvert \int_{\llvert  z'\rrvert  >\beta\llvert  u\rrvert  } \bigl\llvert
\hat{h}\bigl(z'\bigr)\bigr\rrvert \,dz' %
\\[-1pt]
\hspace*{-1pt}&&{}+ \int_{\llvert  z\rrvert  > \beta\llvert  u\rrvert  } \biggl\llvert \frac{\psi'(-\zeta-iz)\psi'(w+iu)}{ (\psi(w+iu) - \psi
(-\zeta-iz) )^2}
\\[-1pt]
\hspace*{-1pt}&&\hspace*{81pt}{}-
\frac{1}{ (w+iu + \zeta+ iz )^2} \biggr\rrvert \int_{\llvert  z'\rrvert  >\llvert  z\rrvert  } \bigl\llvert
\hat{h}\bigl(z'\bigr)\bigr\rrvert \,dz' \,dz %
\\[-1pt]
\hspace*{-1pt}&\leq&\frac{c_1}{\llvert  u\rrvert  } \int_{\RR} \bigl\llvert \hat{h}(z)
\bigr\rrvert \,dz
\\[-1pt]
\hspace*{-1pt}&&{} + c_2 \int_{\llvert  z'\rrvert  >\beta\llvert  u\rrvert  } \bigl\llvert
\hat{h}\bigl(z'\bigr)\bigr\rrvert \,dz' \int
_{\llvert  z\rrvert  >\beta\llvert  u\rrvert  } \biggl(\frac{\llvert  z\rrvert   \llvert  u\rrvert  }{(z^2-u^2)^2 + z^2}
\\[-1pt]
\hspace*{-1pt}&&\hspace*{148pt}{}+ \frac{1}{(w+\zeta)^2 + (z+u)^2}
\biggr) \,dz %
\\[-1pt]
\hspace*{-1pt}&\leq&\frac{c_1}{\llvert  u\rrvert  } \int_{\RR} \bigl\llvert \hat{h}(z)
\bigr\rrvert \,dz
+ c_3 \int_{\llvert  z'\rrvert  >\beta\llvert  u\rrvert  } \bigl\llvert
\hat{h}\bigl(z'\bigr)\bigr\rrvert \,dz' \biggl(\int
_{0}^{\infty}\hspace*{-1pt} \frac{z^2}{(z^2-u^2)^2 + u^2} \,dz + c_4
\biggr) %
\\[-1pt]
\hspace*{-1pt}&\leq&\frac{c_1}{\llvert  u\rrvert  } \int_{\RR} \bigl\llvert \hat{h}(z)
\bigr\rrvert \,dz+ c_5 \int_{\llvert  z'\rrvert  >\beta\llvert  u\rrvert  } \bigl\llvert
\hat{h}\bigl(z'\bigr)\bigr\rrvert \,dz', %
\end{eqnarray*}
where we made use of (\ref{eqfourthintaux}).
Using the above estimate and (\ref{eqtailcond}), we conclude that
$F^w$ belongs to $\mathbb{L}^2(\RR)$. Then the standard properties of
Fourier\vspace*{1pt} transform yield that $F^w$ is a Fourier transform of some
$\tilde{g}^w\in\mathbb{L}^2(\RR)$. Moreover, $\tilde{g}^w$ can be
obtained as the $\mathbb{L}^2(\RR)$-limit of functions $x\mapsto
e^{-w x}g^w_r(x)$, as $r\rightarrow\infty$. Since we showed in step~1
that $g^w_r$ converges weakly to $g$ (which is independent of $w$), we
conclude that $g(x)$ is locally integrable, and $e^{-wx}g(x)$ is square
integrable over $x\in\RR$, for any $w\geq\gamma$. In addition, the
Fourier transform of $x\mapsto e^{-w x}g(x)$ is $F^w$, which implies
\[
\lim_{R\rightarrow\infty}\int_{-R}^R
e^{-(w+iu)x} g(x) \,dx = F^w(u), %
\]
where the convergence is understood in an $\mathbb{L}^2(\RR)$ sense.
Moreover, since $x\mapsto e^{-w x}g(x)$ is square integrable over $x\in
\RR$ and $g$ has support in $[0,\infty)$, it is easy to deduce that
the left-hand side of the above equation converges point wise to a
continuous function of $w+iu$, for all $u\in\RR$ and $w\geq\gamma$.
This, along with Corollary~\ref{coranExt}, implies
\begin{eqnarray*}
\int_{0}^{\infty} e^{-w x} g(x) \,dx &=&
F^w(0)
\\
&=& \int_{\RR} \biggl( \frac{\psi'(w)}{\psi(w) - \psi(-\zeta-iz)}
- \frac{1}{w + \zeta+ iz} \biggr) \hat{h}(z) \,dz
\\
&=& \psi'(w)\Upsilon
\bigl(\psi(w)\bigr),
\end{eqnarray*}
for all $w\geq\gamma$.
As shown in Section~\ref{subseLaplace}, a change of variables
turns the left and the right-hand sides of the above equation into the
Laplace transforms of $\EE g(X_t)$ and $\EE h(X_t)$, respectively. Due
to the uniqueness of the Laplace inverse, the expectations have to
coincide for all $t>0$.

\item[\textit{Step} 4.] It only remains to show that even without the additional
assumptions, $\hat{h}\in\mathbb{L}^1(\RR)$ and (\ref
{eqtailcond}), the action of the generalized function $g\in\mathcal
{D}^*$ (constructed in step~1)\vspace*{1pt} on $p_t$ coincides with $\EE h(X_t)$,
for all $t>0$.
Since $\hat{h}\in\mathbb{L}^1(\RR)$ [resp., $\hat{h}\in\mathbb
{L}^2(\RR)$], there exists a sequence of functions $\hat{h}^n$, such
that every $\hat{h}^n$ is infinitely smooth, with compact support in
$(-\infty,0)$, and $\hat{h}^n$ converges to $\hat{h}$ in $\mathbb
{L}^1(\RR)$ [resp., $\mathbb{L}^2(\RR)$].
Denote by $e^{\zeta x}h^n(x)$ the Fourier transform of $\hat{h}^n$.
Then $e^{\zeta x}h^n(x)$ converges to $e^{\zeta x}h(x)$ in $\mathbb
{L}^{\infty}(\RR)$ [resp., $\mathbb{L}^2(\RR)$].
Notice that every $\hat{h}^n$ belongs to $\mathbb{L}^1(\RR)\cap
\mathbb{L}^2(\RR)$ and satisfies (\ref{eqtailcond}).
Then as shown in steps 2~and~3, there exists a locally integrable
function $g^n$, with at most exponential growth at infinity, which is
the weak symmetry image of $h^n$. Namely, it satisfies
\[
\bigl\langle g^n, p_t \bigr\rangle= \EE
g^n(X_t) = \EE h^n(X_t)\qquad
\forall t>0. %
\]
Notice that, since $X_t$ has a continuous density and since $g^n$ is
locally integrable with at most exponential growth, the expectation of
$g^n(X_t)$ is well defined.
Due to Lemma~\ref{leptl2} and the choice of $h^n$, we obtain
\begin{eqnarray*}
\hspace*{-3pt}&& \EE\bigl\llvert h(X_t) - h^n(X_t)\bigr
\rrvert
\\
\hspace*{-3pt}&&\qquad = \int_{\RR} \bigl\llvert h(x)-h^n(x)
\bigr\rrvert e^{\zeta x} e^{-\zeta x}p_t(x) \,dx %
\\
\hspace*{-3pt}&&\qquad\leq\min \bigl(\bigl\llVert \bigl(h-h^n\bigr) e^{\zeta\cdot}\bigr
\rrVert _{\mathbb{L}^{\infty
}(\RR)} \bigl\llVert p_t e^{-\zeta\cdot}\bigr
\rrVert _{\mathbb{L}^1(\RR)}, \bigl\llVert \bigl(h-h^n\bigr)
e^{\zeta\cdot}\bigr\rrVert _{\mathbb{L}^2(\RR)} \bigl\llVert p_t
e^{-\zeta\cdot}\bigr\rrVert _{\mathbb{L}^2(\RR)} \bigr)
\\
\hspace*{-3pt}&&\qquad\rightarrow0, %
\end{eqnarray*}
as $n\rightarrow\infty$.
Let us show that $\langle g^n, p_t \rangle\rightarrow\langle g, p_t
\rangle$:
\begin{eqnarray*}
&& \bigl\llvert \bigl\langle g - g^n, p_t \bigr\rangle
\bigr\rrvert
\\
&&\qquad \leq
\frac{1}{2\pi} \int_{\RR} \biggl\llvert \int
_{\RR} e^{(\gamma+ iu) y} p_t(y) \,dy\biggr\rrvert
\\
&&\hspace*{60pt}{}\times
\biggl\llvert \int_{\RR} \biggl( \frac{\psi'(\gamma+ iu)}{\psi(\gamma+
iu) - \psi(-\zeta-iz)} -
\frac{1}{\gamma+ iu + \zeta+ iz} \biggr)\biggr\rrvert
\\
&&\hspace*{60pt}{}\times \bigl\llvert \hat{h}(z) -
\hat{h}^n(z)\bigr\rrvert \,dz \,du. %
\end{eqnarray*}
%
Making use of (\ref{eqephiintest}) and (\ref{eqcor2est1}), we
conclude that
\[
\bigl\llvert \bigl\langle g - g^n, p_t \bigr\rangle
\bigr\rrvert \leq c_{1} \bigl(\bigl\llVert \hat{h}-
\hat{h}^n\bigr\rrVert _{\mathbb{L}^1(\RR)} \wedge \bigl\llVert \hat{h}-
\hat{h}^n\bigr\rrVert _{\mathbb{L}^2(\RR)} \bigr)\rightarrow0,
\]
as $n\rightarrow\infty$, which implies $\langle g, p_t \rangle= \EE
h(X_t)$, for all $t>0$, and completes the proof of the theorem.\quad\qed
\end{longlist}\noqed
\end{pf*}

\begin{pf*}{Proof of Corollary~\ref{cormainrate}}
First, we notice that
\[
\bigl\llvert g(x) - g_{r,R}(x) \bigr\rrvert \leq\bigl\llvert g(x) -
g_{r}(x) \bigr\rrvert + \bigl\llvert g_{r}(x) -
g_{r,R}(x) \bigr\rrvert. %
\]
The first term on the right-hand side of the above is bounded by the
right-hand side of (\ref{eqmainrategrR}) due to Theorem~\ref
{thmain}. To analyze the second term, we only need to estimate the
right-hand side of (\ref{eqLevyfinalgRconv}), with the integration
over $z\in\RR$ replaced by the integration over $\llvert  z\rrvert  >R$. We will
refer to it as the modified right-hand side of (\ref
{eqLevyfinalgRconv}). The estimation is done as in step~2 in the
proof of Theorem~\ref{thmain}, with the exception that in the present
case, all the terms vanish, as $r,R\rightarrow\infty$.
Following the derivation of~(\ref{eqfirstintrate}), we conclude
that the first integral on the modified right-hand side of (\ref
{eqLevyfinalgRconv}) is estimated from above by
\begin{eqnarray*}
&& c_1 e^{\gamma x} \int_{\llvert  z\rrvert  >R}
\frac{r}{\sqrt{(z^2 - r^2)^2 + z^2}} \bigl\llvert \hat{h}(z) \bigr\rrvert \,dz
\\
&&\qquad \leq c_1
e^{\gamma x} \int_{R<\llvert  z\rrvert  \leq r/2} \frac{2}{r} \bigl\llvert
\hat{h}(z) \bigr\rrvert \,dz + 2c_1 e^{\gamma x} \int
_{\llvert  z\rrvert  >R,\llvert  z\rrvert  >r/2} \bigl\llvert \hat{h}(z) \bigr\rrvert \,dz. \end{eqnarray*}
Similarly to (\ref{eqsecondintrate}), the second integral on the
modified right-hand side of (\ref{eqLevyfinalgRconv}) is bounded by
\[
2 e^{\gamma x} \int_{R<\llvert  z\rrvert  \leq r/2} \frac{2}{r} \bigl
\llvert \hat{h}(z) \bigr\rrvert \,dz + 2 e^{\gamma x} \int_{\llvert  z\rrvert  >R,\llvert  z\rrvert  > r/2}
\frac{1}{\gamma+ \zeta} \bigl\llvert \hat{h}(z) \bigr\rrvert \,dz. %
\]
Following (\ref{eqthirdintrate}), we estimate the third integral on
the modified right-hand side of~(\ref{eqLevyfinalgRconv}) via
\begin{eqnarray*}
&& c_2 e^{\gamma x} \int_{R<\llvert  z\rrvert  \leq r/2} \int
_{r}^{\infty} \frac{1}{(u-r/2)^2 } \,du \bigl\llvert
\hat{h}(z) \bigr\rrvert \,dz
\\
&&\qquad{}+ c_3 e^{\gamma x} \int
_{\llvert  z\rrvert  >R,\llvert  z\rrvert  > r/2} \int_{0}^{\infty}
\frac{1}{\sqrt{u^4+1}} \,du \bigl\llvert \hat{h}(z) \bigr\rrvert \,dz. %
\end{eqnarray*}
Similarly to (\ref{eqfourthintrate}), we find the upper bound for
the fourth integral on the modified right-hand side of (\ref
{eqLevyfinalgRconv}):
\begin{eqnarray*}
&& c_{4}e^{\gamma x} \int_{R<\llvert  z\rrvert  \leq r/2} \int
_{r}^{\infty} \frac
{u^2}{(u - r/2 )^4} \,du \bigl\llvert
\hat{h}(z) \bigr\rrvert \,dz %
\\
&&\qquad{}+ c_{4}e^{\gamma x} \int_{\llvert  z\rrvert  >R, \llvert  z\rrvert  >r/2} \biggl( 2
(1+1/N )^2 + \int_{\llvert  u\rrvert  >1} \frac{1}{u^2} \,du
\biggr) \bigl\llvert \hat{h}(z) \bigr\rrvert \,dz. %
\end{eqnarray*}
Finally, we obtain the estimate of the fifth integral on the modified
right-hand side of (\ref{eqLevyfinalgRconv}), following the
derivation of (\ref{eqfifthint}):
\begin{eqnarray*}
&& e^{\gamma x} \int_{R<\llvert  z\rrvert  \leq r/2} \int_{r}^{\infty}
\frac
{1}{(u-r/2)^2} \,du \bigl\llvert \hat{h}(z) \bigr\rrvert \,dz
\\
&&\qquad{} + e^{\gamma x}
\int_{\llvert  z\rrvert  >R, \llvert  z\rrvert  > r/2} \int_{\RR} \frac{1}{(\gamma+
\zeta)^2 + u^2}
\,du \bigl\llvert \hat{h}(z) \bigr\rrvert \,dz. %
\end{eqnarray*}
Collecting the above estimates, we obtain (\ref{eqmainrategrR}).
\end{pf*}
\end{appendix}

\section*{Acknowledgments}
We thanks the anonymous referees for their feedback which helped us
improve our paper.



\printaddresses
\end{document}